\documentclass[11pt, leqno, oneside]{amsart} 

\usepackage{amsmath,amsthm,amssymb} 
\usepackage{amsfonts} 
\usepackage{txfonts} 
\usepackage{thmtools} 

\usepackage{geometry} 

\usepackage[dvipsnames, svgnames]{xcolor} 

\usepackage[shortlabels]{enumitem} 

\usepackage[colorlinks,citecolor=blue,hypertexnames=false,urlcolor=blue,breaklinks=true, pagebackref]{hyperref} 

\usepackage{accents} 
\usepackage{exscale} 
\usepackage{dsfont} 
\usepackage{mathrsfs} 
\usepackage[T1]{fontenc} 

\usepackage{etoolbox}
\appto\appendix{\addtocontents{toc}{\protect\setcounter{tocdepth}{1}}}

\renewcommand*\backref[1]{\ifx#1\relax \else (p. #1) \fi} 

\numberwithin{equation}{section} 

\theoremstyle{plain}
\newtheorem{theorem}[equation]{Theorem} 
\newtheorem{lemma}[equation]{Lemma}

\newtheorem{proposition}[equation]{Proposition}

\theoremstyle{definition}
\newtheorem{definition}[equation]{Definition}

\theoremstyle{remark}

\newtheorem{claim*}{Claim}

\newcommand{\bb}[1]{\mathbb{#1}}

\newcommand{\dist}{\operatorname{dist}}

\newcommand{\dv}{\operatorname{div}}
\newcommand{\Di}{\operatorname{D}}
\newcommand{\N}{\operatorname{N}}
\newcommand{\Reg}{\operatorname{R}}
\newcommand{\n}[1]{\mathscr{#1}}
\newcommand{\m}[1]{\mathcal{#1}}
\newcommand{\Lip}{\operatorname{Lip}} 
\newcommand{\ep}{\varepsilon}
\newcommand{\ra}{\rightarrow}
\newcommand{\lra}{\longrightarrow}   
\newcommand{\loc}{\operatorname{loc}}

\newcommand{\wt}{\widetilde}

\newcommand{\R}{\mathbb{R}}

\def\R{\mathbb{R}}
\def\pom{{\partial\Omega}}

\def\ve{\varepsilon}

\DeclareMathOperator{\supp}{supp}
\DeclareMathOperator{\diam}{diam}

\def\Xint#1{\mathchoice
	{\XXint\displaystyle\textstyle{#1}}%
	{\XXint\textstyle\scriptstyle{#1}}%
	{\XXint\scriptstyle\scriptscriptstyle{#1}}%
	{\XXint\scriptscriptstyle\scriptscriptstyle{#1}}%
	\!\int}
\def\XXint#1#2#3{{\setbox0=\hbox{$#1{#2#3}{\int}$}
		\vcenter{\hbox{$#2#3$}}\kern-.5\wd0}}

\def\dashint{\Xint-}

\def\YYint#1#2#3{{\setbox0=\hbox{$#1{#2#3}{\iint}$}
		\vcenter{\hbox{$#2#3$}}\kern-.51\wd0}}
 
\newcommand*{\dt}[1]{%
	\accentset{\mbox{\Large\bfseries .}}{#1}}

\allowdisplaybreaks

\setlist{nosep}

\begin{document}
	
\author[M. Mourgoglou]{Mihalis Mourgoglou}
\email[]{michail.mourgoglou@ehu.eus}
\address{Departamento de Matemáticas, Universidad del País Vasco, UPV/EHU, Barrio Sarriena S/N 48940 Leioa, Spain and, IKERBASQUE, Basque Foundation for Science, Bilbao, Spain}

\author[B. Poggi]{Bruno Poggi}
\email[]{brunopoggi@pitt.edu}
\address{Department of Mathematics, University of Pittsburgh, Pittsburgh, PA, 15213}

 \author[M. Ulmer]{Martin Ulmer}
\email[]{martin\_ulmer@brown.edu}
\address{Department of Mathematics, Brown University, 151 Thayer Street, Providence, RI, 02912, USA}

\title[On boundaries of boundary value problems]{On boundaries of boundary value problems: The cases of Regularity and Neumann}
\date{\today}

\thanks{M.M. was supported  by IKERBASQUE and partially supported by the grant PID2024-157724NB-I00: of the Ministerio de Econom\'ia y Competitividad (Spain), and by  IT-1615-22 (Basque Government). BP was partially supported by an award from the Simons Foundation in the Travel Support for Mathematicians program, [SFI-MPS-TSM-00025559, BP]. 
}
 
\begin{abstract} We show that on a bounded uniform domain $\Omega\subset\bb R^{n+1}$, $n\geq2$, with $n$-Ahlfors regular boundary, and given a second-order homogeneous linear strongly elliptic operator $L=-\dv A\nabla$ with bounded measurable coefficients, for which both  $(\Di_{p'}^{L^*})$ and $(\N_p^L)$ are solvable in $\Omega$, the following sharp and novel approximation result holds for the  Neumann problem 
\[
\left\{\begin{aligned}-\dv A\nabla u&=0,\qquad&\text{in }&\Omega,\\ \partial_{\nu_A} u&=g,\qquad&\text{on }&\partial\Omega,\end{aligned}\right. 
\]
with data $g\in L^p$ $p\in(1,\infty)$: There is a sequence of compactly supported Lipschitz functions $\{H_j\}_j$ with uniformly bounded norm in the space ${\bf C}_{2,p}$, such that if   $w_j$ is the weak solution $w_j\in W^{1,2}(\Omega)$ to the Poisson-Neumann problem
\[
\left\{\begin{aligned}-\dv A\nabla w_j&= H_j,\qquad&\text{in }&\Omega,\\ \partial_{\nu_A}w_j&=0,\qquad&\text{on }&\partial\Omega,\end{aligned}\right. 
\]
then $\nabla w_j\ra \nabla u$ as $j\ra\infty$ in the weak-$*$ topology of $({\bf N}_{2,p})^{n+1}$. We also show a variant of this approximation theorem for the rough Neumann problem with data in $W^{-1,p}$. These results are   analogues of an approximation result recently discovered by the first two authors in the case of the Dirichlet problem with data in $L^p$. 

Moreover, we show that an analogous approximation theorem is not viable for the Dirichlet-Regularity problem, nor for the variational Dirichlet problem. All of the results are new even in the case of the Laplacian and on the unit ball.

\end{abstract}

\maketitle

{\bf Statement on human authorship and use of artificial intelligence.} All novel results, arguments, and exposition in this article are due  to the intellectual contributions of the human authors. Artificial intelligence tools were used   solely for minor linguistic editing and proofreading at the stage of the final draft.

\hypersetup{linkcolor=blue}

	\tableofcontents 

\section{Introduction}

This manuscript develops a program started in \cite{mp26}, which presented a novel approximation theorem for the homogeneous Dirichlet problem for elliptic PDEs with data in $L^p$ by inhomogeneous Poisson problems with trivial boundary data (see Theorem \ref{thm.approx} recapped below). The results in \cite{mp26} left open the question as to whether the approximation result has analogues in simpler settings, for instance that of the  variational Dirichlet problem. Another natural line of research is to seek analogues of Theorem \ref{thm.approx} for other boundary value problems. Here, we study such analogues   for the Dirichlet problem at higher orders of regularity (in particular, the so-called Regularity problem) and for the Neumann problem. The results we find present a story which is not straightforward.

\subsection{Motivation.} But let us begin with the original motivation for seeking such approximation theorems. The $L^p$-solvability of homogeneous boundary value problems for second-order linear  elliptic partial differential equations   with non-smooth coefficients  on rough domains  has been an area of intense and fruitful research in the last few decades, culminating in shocking equivalences between certain geometric properties of domains   in the Euclidean space  and $L^p$-solvability of the homogeneous Dirichlet problem   for the Laplacian. In \cite{mpt25}, the authors discovered under no  assumptions on the real elliptic operator  and quite relaxed conditions on the geometry of the domain that the solvability of the Dirichlet problem with data in $L^p$
\begin{equation}\label{eq.dirichlet}
(\operatorname{D}^L_p)\qquad\left\{\begin{aligned}-\dv A\nabla u&=0,\quad&\text{in }&\Omega,\\u&=g,\quad&\text{on }&\partial\Omega,\\ \Vert\m N(u)\Vert_{L^p(\partial\Omega,\sigma)}&\lesssim\Vert g\Vert_{L^p(\partial\Omega,\sigma)},&&\end{aligned}\right.
\end{equation}
is equivalent to solvability of the Poisson-Dirichlet problem
\begin{equation}\label{eq.poisson}
(\operatorname{PD}^L_p)\qquad\left\{\begin{aligned}-\dv A\nabla w&=-\dv F,\qquad&\text{in }&\Omega,\\ w&=0,\qquad&\text{on }&\partial\Omega,\\  \Vert\widetilde{\m N}_2(w)\Vert_{L^p(\partial\Omega,\sigma)}&\lesssim\Vert\n C_2(F)\Vert_{L^p(\partial\Omega,\sigma)},&&\end{aligned}\right. 
\end{equation}
with inhomogeneous data $F$ in the  $L^p$-Carleson space ${\bf C}_{2,p}$. For the definitions of the non-tangential maximal operators $\m N$ and $\widetilde{\m N}_2$, as well as that of the Carleson function $\n C_2$ and the space ${\bf C}_{2,p}$, see Section \ref{sec.ntmax} and Definition \ref{def.spaces}. This quantitative equivalence of the problems (\ref{eq.dirichlet}) and (\ref{eq.poisson}) has garnered interest in the mathematical community and is an ongoing area of research across different settings and boundary value problems \cite{mz25,fl24,bmp25,fyy25}.

In the recent manuscript \cite{mp26}, the authors investigated more deeply the relationship between the boundary value problems (\ref{eq.dirichlet}) and (\ref{eq.poisson}), and found the following surprising approximation result. We defer all definitions to Section \ref{sec.prelim}.

\begin{theorem}[Approximation Theorem for the Dirichlet problem \cite{mp26}]\label{thm.approx} Let $\Omega\subset\bb R^{n+1}$, $n\geq2$, be a bounded   domain satisfying the corkscrew condition and with $n$-Ahlfors regular boundary, and let $A$ be a strongly elliptic, not necessarily symmetric real matrix of bounded measurable coefficients in $\Omega$. Fix $p\in(1,\infty)$, and suppose that $(\Di_p^L)$ is solvable in $\Omega$. Fix $g\in L^p(\partial\Omega)$, and let $u\in{\bf N}_{2,p}\cap W^{1,2}_{\loc}(\Omega)$ be the unique  solution to the   Dirichlet problem 
	\begin{equation}\label{eq.dirichlet0}
		\left\{\begin{aligned}-\dv A\nabla u&=0,\qquad&\text{in }&\Omega,\\ u&=g,\qquad&\text{on }&\partial\Omega.\end{aligned}\right. 
	\end{equation}
	Then  there exists a sequence of vector functions $\{F_j\}_{j\in\bb N}\subset\Lip_c(\Omega)^{n+1}$, with the following properties:
	\begin{enumerate}[(i)]
		\item\label{item.c1} There exists $C>0$ depending only on dimension, ellipticity, the corkscrew   constant, and the $n$-Ahlfors regularity constant, such that
		\begin{equation}\label{eq.uniformc}
			\Vert F_j\Vert_{{\bf C}_{2,p}}\leq C\Vert g\Vert_{L^p(\partial\Omega,\sigma)}.
		\end{equation}
		\item\label{item.c2} For each $j\in\bb N$, let $w_j$ be the unique weak solution in $W_0^{1,2}(\Omega)$ to the Poisson problem
		\begin{equation}\label{eq.poissondirichlet0}
			\left\{\begin{aligned}-\dv A\nabla w_j&=-\dv F_j,\qquad&\text{in }&\Omega,\\ w_j&=0,\qquad&\text{on }&\partial\Omega.\end{aligned}\right. 
		\end{equation}
		Then, $w_j\ra u$ as $j\ra\infty$ in the following topologies:
		\begin{enumerate}[(a)]
			\item\label{item.c2a} the strong topology in $C_{\loc}^{\alpha}(\Omega)$ for some $\alpha\in(0,1)$.
			\item\label{item.c2b} the strong topology in $W^{1,2}_{\loc}(\Omega)$.
			\item\label{item.c2c} the weak-$*$ topology in ${\bf N}_{2,p}$. Moreover, this convergence cannot be improved to convergence in the weak topology in ${\bf N}_{2,p}$ unless $g\equiv0$.
		\end{enumerate}
	\end{enumerate}
\end{theorem}

This theorem may be interpreted as saying that the solutions to the Dirichlet problem (\ref{eq.dirichlet}) lie on the boundary (in the weak-$*$ topology of ${\bf N}_{2,p}$) of the solution space to the Poisson-Dirichlet problem (\ref{eq.poisson}) (see \cite[Corollary 1.12]{mp26}). See the Introduction of \cite{mp26} for an extensive discussion on this result, including why the result is not feasible for the Dirichlet problem (\ref{eq.dirichlet0}) with continuous boundary data.

\subsection{Main results.} It is natural to wonder whether an analogue of Theorem \ref{thm.approx} holds at different levels of regularity of the data of the Dirichlet problem. One such instance is the well-known Dirichlet Regularity problem,
\begin{equation}\label{eq.regularity}
(\operatorname{R}^L_p)\qquad\left\{\begin{aligned}-\dv A\nabla u&=0,\quad&\text{in }&\Omega,\\u&=f,\quad&\text{on }&\partial\Omega,\\ \Vert\m N(\nabla u)\Vert_{L^p(\partial\Omega,\sigma)}&\lesssim\Vert f\Vert_{\dt W^{1,p}(\partial\Omega,\sigma)},&&\end{aligned}\right.
\end{equation}
where $\dt W^{1,p}(\partial\Omega,\sigma)$ is a Haj\l{}asz Sobolev space on $\partial\Omega$. For the definition of Haj\l{}asz Sobolev spaces,  see Section \ref{sec.sobolev}; for more details on the definition of the Regularity problem, see Definition \ref{def.regularity}. Another important and classical level of regularity is that of the variational Dirichlet problem, whereby one seeks to find a weak solution $u\in W^{1,2}(\Omega)$ to
\begin{equation}\label{eq.variational}
		(\operatorname{D}^L_{\operatorname{var}})\qquad\left\{\begin{aligned}-\dv A\nabla u&=0,\qquad&\text{in }&\Omega,\\ u&=f,\qquad&\text{on }&\partial\Omega,\end{aligned}\right. 
\end{equation}
for\footnote{Here, we are assuming $\Omega$ is a bounded uniform domain with $n$-Ahlfors regular boundary, so that the identification of the fractional Sobolev space $H^{\frac12}(\partial\Omega)$ as the trace space of $W^{1,2}(\Omega)$ is not an issue.} $f\in H^{\frac12}(\partial\Omega)$. Indeed, the solvability of the variational Dirichlet problem in this setting is a classical and basic result in elliptic PDE theory, and there is also an extensive literature on the Dirichlet Regularity problem \cite{mt24, mpt25, dhp22}  

One might hope that Theorem \ref{thm.approx} attains an appropriate analogue for one or both of the problems in (\ref{eq.regularity}) or (\ref{eq.variational}). However, it is remarkable that, despite Theorem \ref{thm.approx}, neither the Regularity problem nor the variational Dirichlet problem admit a sharp approximation by solutions to the Poisson problem with trivial boundary data, essentially due to the same obstruction.

\begin{proposition}[Failure of Approximation Theorem for the Regularity Problem]\label{prop.approxReg} Let $\Omega\subset\bb R^{n+1}$, $n\geq2$, be a bounded   domain satisfying the corkscrew condition and with uniformly $n$-rectifiable boundary, and let $A$ be a strongly elliptic, not necessarily symmetric real matrix of bounded measurable coefficients in $\Omega$. Fix $p\in(1,\infty)$ and suppose that $(\operatorname{R}^L_p)$ is solvable in $\Omega$. Let $f\in\dt W^{1,p}(\partial\Omega)$, and let $u\in W^{1,2}_{\loc}(\Omega)$  be the unique  solution to the Regularity problem (\ref{eq.regularity}).

	  If there exists a sequence of functions  $\{w_j\}_j\subset W_0^{1,2}(\Omega)$ such that
     $\nabla w_j\lra \nabla u$ as $j\ra\infty$ in the weak-$*$ topology of $({\bf N}_{2,p})^{n+1}$, then $f$ is a constant function.
\end{proposition}

\begin{proposition}[Failure of Approximation Theorem for the variational Dirichlet Problem]\label{prop.approxvar} Let $\Omega\subset\bb R^{n+1}$, $n\geq2$, be a bounded  uniform domain with $n$-Ahlfors regular boundary, and let $A$ be a strongly elliptic, not necessarily symmetric real matrix of bounded measurable coefficients in $\Omega$.   Let $f\in H^{\frac12}(\partial\Omega)$, and let $u\in W^{1,2}(\Omega)$   be the unique  solution to the variational Dirichlet problem (\ref{eq.variational}).

	  If there exists a sequence of functions  $\{w_j\}_j\subset W_0^{1,2}(\Omega)$ such that
     $\nabla w_j\lra \nabla u$ as $j\ra\infty$ in the weak\footnote{Note that the weak-$*$ topology and the weak topology of $L^2(\Omega)$ coincide.}  topology of $(L^2(\Omega))^{n+1}$, then $f$ is a constant function.
\end{proposition}
 
 These two results show that Theorem \ref{thm.approx}  is necessarily a low-order phenomenon and a distinguishing feature of the $L^p$ Dirichlet problem. Since an analogue of Theorem \ref{thm.approx} is also not feasible for the Dirichlet problem with continuous boundary data (see the discussion in the introduction of \cite{mp26}), these arguments explain why such an approximation theorem was not revealed in the literature in more basic settings before; indeed, it is simply \emph{not true} in such settings. The proofs of Propositions \ref{prop.approxReg} and  \ref{prop.approxvar} are short and  similar; they are given in Section \ref{sec.failure}.

Now let us turn to positive results in a different direction. Theorem \ref{thm.approx} says, roughly speaking, that the $L^p$ Dirichlet boundary value problem is ``on the boundary'' of the Poisson (Dirichlet) boundary value problem with Carleson $L^p$ data. It is a natural question to consider whether this structure survives to other classical boundary value problems, such as the Neumann problem. 

In light of the negative results of Propositions \ref{prop.approxReg} and \ref{prop.approxvar}, it is not necessarily a priori clear that there is an approximation theorem for the $L^p$ Neumann problem, since its data lives at the same level of regularity as the Dirichlet Regularity problem. On the other hand, since the solutions to the Poisson problem with trivial Neumann boundary data do not necessarily have null trace on the boundary, the main obstacle exploited in the proof of Proposition  \ref{prop.approxReg} is not available in the setting of the Neumann problem. 

The second set of main results in this manuscript shows that, indeed, this structure unveiled for the $L^p$ Dirichlet problem in Theorem \ref{thm.approx}  survives to the case of the Neumann problem. In fact, here we are able to deduce approximation theorems at two separate levels of regularity of the data. In preparation for the statements of the approximation theorems for the Neumann problem, let us denote by $\widehat W^{1,2}(\Omega)$  the set of functions in $W^{1,2}(\Omega)$ with zero mean trace in $\partial\Omega$.  

The theorem below describes the approximability of solutions to the $L^p$ Neumann problem by solutions to Poisson-Neumann problems with trivial Neumann boundary data, in analogy to Theorem \ref{thm.approx}.

\begin{theorem}[Approximation Theorem for the $L^p$ Neumann problem]\label{thm.approxNeumann} Let $\Omega\subset\bb R^{n+1}$, $n\geq2$, be a bounded uniform domain with $n$-Ahlfors regular boundary, and let $A$ be a strongly elliptic, not necessarily symmetric real matrix of bounded measurable coefficients in $\Omega$. Fix $p\in(1,\infty)$, and suppose that $(\Di_{p'}^{L^*})$ and $(\N_p^L)$ are solvable in $\Omega$. Fix $g\in L^p(\partial\Omega)$ with $\int_{\partial\Omega} g\,d\sigma=0$, and let $u\in W^{1,2}_{\loc}(\Omega)$ with $\nabla u\in ({\bf N}_{2,p})^{n+1}$ be the unique  solution to the Neumann problem 
	\begin{equation}\label{eq.neumann0}
		\left\{\begin{aligned}-\dv A\nabla u&=0,\qquad&\text{in }&\Omega,\\ \partial_{\nu_A} u&=g,\qquad&\text{on }&\partial\Omega.\end{aligned}\right. 
	\end{equation}
	Then  there exists a sequence of functions $\{H_j\}_{j\in\bb N}\subset\Lip_c(\Omega)$, with the following properties:
	\begin{enumerate}[(i)]
		\item\label{item.c1} There exists $C>0$ depending only on dimension, ellipticity, the corkscrew   constant, and the $n$-Ahlfors regularity constant, such that
		\begin{equation}\label{eq.uniformcNeu}
			\Vert H_j\Vert_{{\bf C}_{2,p}}\leq C\Vert g\Vert_{L^p(\partial\Omega,\sigma)}.
		\end{equation}
		\item\label{item.c2} For each $j\in\bb N$, let $w_j$ be the unique weak solution $w_j\in\widehat{W}^{1,2}(\Omega)$ to the Poisson-Neumann problem
		\begin{equation}\label{eq.poissonNeumann0}
			\left\{\begin{aligned}-\dv A\nabla w_j&= H_j,\qquad&\text{in }&\Omega,\\ \partial_{\nu_A}w_j&=0,\qquad&\text{on }&\partial\Omega.\end{aligned}\right. 
		\end{equation}
		Then, $\nabla w_j\ra \nabla u$ as $j\ra\infty$ in the weak-$*$ topology in $({\bf N}_{2,p})^{n+1}$. In particular, $\nabla w_j\ra \nabla u$ as $j\ra\infty$ strongly in $L_{\loc}^2(\Omega)$.
	\end{enumerate}
\end{theorem}

The proof of Theorem \ref{thm.approxNeumann} is given in Section \ref{sec.neumann}. The $L^{p}$ Neumann problem for bounded uniform domains with $n$-uniformly rectifiable boundary remains unsolved even for the Laplacian, although there are results in Lipschitz domains in \cite{jk81} and certain rough domains \cite{mt26}.  A ``dual'' problem to the $L^p$ Neumann problem has been identified   as the ``rough Neumann problem'' with data in a ``negative'' Sobolev space (see \cite{AM14,mt26}). This problem has been useful in establishing the solvability of the $L^p$ Neumann problem in certain settings \cite{mt26}, in much the same way as the $L^p$ Dirichlet problem is useful in solving the $L^{p'}$ Regularity problem.  Our final main result concerns the approximability of the solutions to the rough Neumann problem by solutions to the Poisson-Neumann problem with trivial Neumann boundary data. We state it in high geometric generality; for the definitions, including the definition of the measure theoretic boundary $\partial_*\Omega$, see Section \ref{sec.geom}.

\begin{theorem}[Approximation Theorem for the rough Neumann problem]\label{thm.approxRoughNeumann} Let $\Omega\subset\bb R^{n+1}$, $n\geq2$, be a bounded uniform domain with $n$-rectifiable and $n$-Ahlfors regular boundary. Let $A$ be a strongly elliptic, not necessarily symmetric real matrix of bounded measurable coefficients in $\Omega$. Fix $p\in(1,\infty)$, and suppose that $(\Di_p^L)$ and $(\N_{p'}^{L^*})$ are solvable in $\Omega$. Fix $g\in L^p(\partial^*\Omega)$ and $1\leq l,k\leq n+1$, and let $u\in W^{1,2}_{\loc}(\Omega)$ with $u\in {\bf N}_{2,p}$ and $\int_\Omega u\, dm=0$ be the unique  solution to the rough Neumann problem 
	\begin{equation}\label{eq.neumann0Rough}
		\left\{\begin{aligned}-\dv A\nabla u&=0,\qquad&\text{in }&\Omega,\\ \partial_{\nu_A} u&=\partial_{t,l,k}g,\qquad&\text{on }&\partial^*\Omega.\end{aligned}\right. 
	\end{equation}
	Then there exists two sequences of functions $\{H_{j}^{(k)},H_j^{(l)}\}_{j\in\bb N}\subset\Lip_c(\Omega)$, with the following properties:
	\begin{enumerate}[(i)]
		\item\label{item.c1} There exists $C>0$ depending only on dimension, ellipticity, the corkscrew   constant, and the $n$-Ahlfors regularity constant, such that
		\begin{equation}\label{eq.uniformc}
			\Vert H_j^{(l)}\Vert_{{\bf C}_{2,p'}} + \Vert H_j^{(k)}\Vert_{{\bf C}_{2,p'}}\leq C\Vert g\Vert_{L^p(\partial\Omega,\sigma)}.
		\end{equation}
		\item\label{item.c2} For each $j\in\bb N$, let $w_j$ be the unique weak solution in $\widehat{W}^{1,2}(\Omega)$ to the Poisson problem
		\begin{equation}\label{eq.poissonNeumann0Rough}
			\left\{\begin{aligned}-\dv A\nabla w_j&= \partial_k H_j^{(k)}-\partial_l H_j^{(l)},\qquad&\text{in }&\Omega,\\ \partial_{\nu_A}w_j&=0,\qquad&\text{on }&\partial\Omega.\end{aligned}\right. 
		\end{equation}
		Then, $w_j\ra u$ as $j\ra\infty$ in the following topologies:
		\begin{enumerate}[(a)]
			\item\label{item.c2a} the strong topology in $C_{\loc}^{\alpha}(\Omega)$ for some $\alpha\in(0,1)$.
			\item\label{item.c2b} the strong topology in $W^{1,2}_{\loc}(\Omega)$.
			\item\label{item.c2c} the weak-$*$ topology in ${\bf N}_{2,p}$.
		\end{enumerate}
	\end{enumerate}
\end{theorem}

We prove Theorem \ref{thm.approxRoughNeumann} in Section \ref{sec.rough}. Let us point out that, in contrast to the negative result for the Dirichlet Regularity problem in Proposition \ref{prop.approxReg}, for the case of the Neumann problem we have the solvability at the classical level of regularity of order 1, as well as the order 0 regularity present in the rough Neumann problem.

The proofs of Theorem~\ref{thm.approxNeumann} and Theorem~\ref{thm.approxRoughNeumann} are inspired by the proof of Theorem~\ref{thm.approx} in \cite{mp26}. Indeed, since the proof of Theorem~\ref{thm.approx} does not rely on the harmonic measure or quantitative estimates for the Green's function, which are tools that are often used for the Dirichlet problem but not available for the Neumann problem, the core integration by parts idea can be adapted to this Neumann setting.

As in \cite{mp26}, we use Varopoulos extensions \cite{var77,var78,hr18,HT21} of the boundary data $g$ to construct explicit interior Poisson data $H_j$, or $H_j^{(k)}$ and $H_j^{(l)}$, for the approximation sequence $w_j$. To obtain a high degree of geometric generality, we rely on the state-of-the-art technology developed in \cite{mz25} (see also Lemma~\ref{lm.varopoulos}). The main additional difficulties come from the variational formulation of the Neumann and rough Neumann problems. In contrast to the Dirichlet problem, where the boundary data are defined pointwise and pointwise convergence of the Varopoulos extension to $g$ is therefore sufficient, the Neumann boundary data enter only through the variational formulation. For the Neumann approximations, we consequently need an extension that converges weakly in $L^p$ on thickenings of the boundary to the prescribed Neumann boundary data.

\subsection{Outline} The rest of this manuscript is organized as follows. In Section \ref{sec.prelim} we introduce all relevant definitions and state useful technical lemmas. In Section \ref{sec.failure} we give the proofs of Propositions \ref{prop.approxReg} and \ref{prop.approxvar}. In Section \ref{sec.neumann} we give the proof of Theorem \ref{thm.approxNeumann}, and in Section \ref{sec.rough} we prove Theorem \ref{thm.approxRoughNeumann}.

\section{Preliminaries}\label{sec.prelim}

\subsection{Geometric preliminaries}\label{sec.geom}

If $x\in\bb R^{n+1}$ and $r>0$, by $B(x,r)$ we denote the Euclidean \emph{open} ball of radius $r$ centered at $x$. If $E\subset\bb R^{n+1}$ and $r>0$, we set $B(E,r):=\{x\in\bb R^{n+1}:\dist(x,E)<r\}$, and ${\bf 1}_E$ is the indicator function of $E$. We  write $a\lesssim b$ if there exists a constant $C>0$ so that $a\leq Cb$ and $a\lesssim_tb$ if $C$ depends on the parameter $t$. We write $a\approx b$ to mean $a\lesssim b\lesssim a$ and define $a\approx_tb$ similarly. Throughout, the letter $C$ denotes a constant that does not depend on the salient parameters in a given proof, and the value of $C$ might change from one line to the next.

We denote by $m=m_{d}$ the Lebesgue measure on $\bb R^{d}$, $d\in\bb N$. 

A set $E\subset\bb R^{n+1}$ is called \emph{$n$-Ahlfors regular} if there exists some constant $C_0>0$ such that
\[
C_0^{-1}r^n\leq\m H^n(B(x,r)\cap E)\leq C_0r^n,\qquad\text{for all }x\in E\text{ and all }0<r\leq\diam(E),
\]
where $\m H^n$ is the $n$-dimensional Hausdorff measure, which we assume to be normalized so that it coincides with $m_n$ in $\bb R^n$.

A set $E\subset\bb R^{n+1}$ is \emph{uniformly $n$-rectifiable} if it is $n$-Ahlfors regular and there exist constants $\theta, M>0$ such that for all $x\in E$ and all $r\in(0,\diam E)$, there is a Lipschitz mapping $g$ from the ball $B_n(0,r)$ in $\bb R^n$ to $\bb R^{n+1}$ with $\operatorname{Lip}(g)\leq M$ such that $\m H^n\big(B(x,r)\cap g(B_n(0,r))\big)\geq\theta r^n$.

By a domain we mean a connected open set. If $\Omega$ is a domain in $\bb R^{n+1}$, we denote by $\sigma=\m H^n|_{\partial\Omega}$ the $n$-dimensional Hausdorff measure restricted to $\partial\Omega$. 

We say that an open set $\Omega$ satisfies the \emph{corkscrew condition} if there exists   $c>0$ such that for each $x\in\partial\Omega$ and every $r\in(0,2\diam\Omega)$, there exists a ball $B\subset B(x,r)\cap\Omega$ so that $r(B)\geq cr$. 

Given two points $x,x'\in\Omega$, and a pair of numbers $M,N\geq1$, an $(M,N)$-\emph{Harnack chain connecting $x$ to $x'$}, is a chain of open balls $B_1,\ldots,B_N\subset\Omega$, with $x\in B_1,x'\in B_N, B_k\cap B_{k+1}\neq\varnothing$ and $M^{-1}\diam B_k\leq\dist(B_k,\partial\Omega)\leq M\diam B_k$. We say that $\Omega$ satisfies the \emph{Harnack chain condition} if there is a uniform constant $M$ such that for any $x,x'\in\Omega$, there exists an $(M,N)$-Harnack chain connecting them, with $N$ depending only on $M$ and on $|x-x'|/\min\{\delta(x),\delta(x')\}$.

We say that a domain $\Omega$ is \emph{uniform} if it satisfies the Harnack chain condition and the corkscrew condition. We say that a domain $\Omega$ is NTA (\emph{non-tangentially accessible}) if it is uniform and $\Omega_{\operatorname{ext}}:=\bb R^{n+1}\setminus\overline\Omega$ satisfies the corkscrew condition. We say a connected open set $\Omega$ is a \emph{chord-arc domain} if it is non-tangentially accessible and $\partial\Omega$ is $n$-Ahlfors regular.

Throughout the rest of this section, unless explicitly mentioned otherwise, $\Omega$ denotes a domain in $\bb R^{n+1}$ satisfying the (interior) corkscrew condition and with $n$-Ahlfors regular boundary.

We now introduce the Hofmann-Mitrea-Taylor tangential derivatives and the required setup, following \cite{hmt10}. First, given a unit vector $\nu\in\bb R^{n+1}$ and a point $x\in\partial\Omega$, set the half-spaces
\[
H^{\pm}_{\nu}(x):=\big\{y\in\bb R^{n+1}~:~\pm\nu\cdot(y-x)\geq0\big\}.
\]
Moreover, set $\Omega^+:=\Omega$, and $\Omega^-:=\bb R^{n+1}\setminus\Omega$. We say that a unit vector $\nu$ in $\bb R^{n+1}$ is the \emph{measure-theoretic outer normal at $\xi$} if 
\[
\lim_{r\ra0}\frac1{r^{n+1}}m\big(B(\xi,r)\cap\Omega^{\pm}\cap H^{\pm}_{\nu}(\xi)\big)=0.
\]  
The measure-theoretic outer normal is unique if it exists. The \emph{measure-theoretic boundary of $\Omega$} is
\begin{equation}\label{eq.measuretheoreticbdd}
\partial_*\Omega:=\Big\{\xi\in\partial\Omega~:~\limsup_{r\ra0}\frac1{r^{n+1}}m\big(B(\xi,r)\cap\Omega^{\pm}\big)>0\Big\}.
\end{equation}
We have that $\sigma$-almost every $\xi\in\partial_*\Omega$ has a measure-theoretic outer normal \cite[p.208]{eg1}, which we denote by $\nu(\xi)=(\nu_j(\xi))_{j=1}^{n+1}$. Moreover, for $\sigma$-almost every $\xi\in\partial_*\Omega$   it also holds true that $|\nu(\xi)|=1$ and that
\[
\nu(\xi)=\lim_{r\ra0}\frac1{\sigma(B(\xi,r))}\int_{B(\xi,r)}\nu(\zeta)\,d\sigma(\zeta).
\]
The \emph{reduced boundary} is given by
\begin{equation}\label{eq.reducedbdd}\partial^*\Omega:=\left\{\xi\in\partial\Omega:\text{there exists } \nu(\xi)\in\mathbb S^n\text{ such that } \lim_{r\to0}\frac{m(B(\xi,r)\cap\Omega^\pm\cap H^\pm_{\nu(\xi)}(\xi))}{r^{n+1}}=0\right\}.
\end{equation}
It follows that \(\partial^*\Omega\subset\partial_*\Omega\), and since \(\Omega\) has finite perimeter, both \(\partial^*\Omega\) and \(\partial_*\Omega\) are countably \(n\)-rectifiable.

Now, given a real-valued function $\varphi\in C_c^1(\bb R^{n+1})$, we define the \emph{Hofmann-Mitrea-Taylor} tangential derivatives of $\varphi$ as
\begin{equation}\label{eq.hmtderivative}
\partial_{t,j,k}\varphi(\xi):=\nu_j(\xi)\partial_k\varphi(\xi)-\nu_k(\xi)\partial_j\varphi(\xi),\qquad x\in\partial_*\Omega,\qquad 1\leq j,k\leq n+1.
\end{equation}
See Sections 3 and 4 of \cite{hmt10} for more on the Hofmann-Mitrea-Taylor tangential derivatives.








Recall that $C_c^{\infty}(\Omega)$ is the space of compactly supported smooth functions in $\Omega$,  and  that for $p\in[1,\infty)$, $W^{1,p}(\Omega)$ is the Sobolev space of $p$-th integrable functions in $\Omega$ whose weak derivatives exist in $\Omega$ and are $p$-th integrable functions, while $W_0^{1,p}(\Omega)$ is the completion of $C_c^{\infty}(\Omega)$ under the norm $\Vert u\Vert_{W^{1,p}(\Omega)}:=\Vert u\Vert_{L^p(\Omega)}+\Vert\nabla u\Vert_{L^p(\Omega)}$. Moreover,   $\dt W^{1,p}(\Omega)$ consists of the $L^1_{\loc}(\Omega)$ functions whose weak gradient is $p$-th integrable over $\Omega$. 

We denote by $\widehat{W}^{1,2}(\Omega)$  the space of all functions in $W^{1,2}(\Omega)$ with zero mean value of the trace on $\partial\Omega$.


If $X$ is a metric space, by $\Lip_c(X)$ we denote the space of Lipschitz functions compactly supported in $X$.


\subsection{Covering $\partial^*\Omega$ by boundaries of Lipschitz subdomains}

In Section \ref{sec.rough}, we shall use a covering of the reduced
boundary $\partial^*\Omega$, up to a set of $\mathcal H^n$-measure
zero, by boundaries of Lipschitz subdomains of $\Omega$. Such a
covering can be obtained by adapting the constructions in
\cite[Theorem 2.5 and Appendix A]{ABHM19} and
\cite[Theorem 1.5]{M21}. More precisely, one first restricts to the
$\sigma$-density points of $\partial^*\Omega$. By the upper Ahlfors
regularity bound,
\[
\sigma\bigl(B(\xi,r)\setminus\partial^*\Omega\bigr)=o(r^n)
\quad\text{as }r\to0
\]
for $\mathcal H^n$-almost every $\xi\in\partial^*\Omega$.
Consequently, at almost every such point, the approximate tangent
plane to $\partial^*\Omega$ is also an approximate tangent plane
to $\partial\Omega$, and the same cone construction applies.
By passing to compact pieces and leaving a margin between the
cone apertures used in the construction, one also obtains the
uniform cone containment property stated below. We record the
properties of this covering that we shall require.

\begin{theorem}\label{thm.corona}
Suppose that $\Omega\subset\bb R^{n+1}$ is a bounded domain with
$n$-Ahlfors regular boundary, and fix $\theta\in(0,\pi/2)$.
Then there exists a countable collection of starlike Lipschitz
domains $\Omega_i\subset\Omega$, with uniform Lipschitz character
after rescaling by their diameters, such that
\begin{equation}\label{eq.exhaustion.reducedboundary}
\mathcal H^n\left(
\partial^*\Omega\setminus\bigcup_i\partial\Omega_i
\right)=0.
\end{equation}
Moreover, there exist a constant $c>0$, independent of $i$, and
unit vectors $e_i\in\bb R^{n+1}$ such that
\begin{equation}\label{eq.exhaustion.cones}
\left\{
\xi+v:
0<|v|<c\,\operatorname{diam}(\Omega_i),\
v\cdot e_i>\cos\theta\,|v|
\right\}
\subset\Omega_i
\end{equation}
for $\mathcal H^n$-almost every
$\xi\in\partial\Omega_i\cap\partial^*\Omega$.
The constant $c$ and the Lipschitz character are controlled
uniformly in terms of $n$, $\theta$, and the Ahlfors regularity
constant of $\partial\Omega$.
\end{theorem}


The following   lemma says that a Lipschitz subdomain of a larger domain shares the outer normal vector with the larger domain on the common boundary. We omit its simple proof.
\begin{lemma}\label{lm.samenormal} Suppose $\Omega\subset\bb R^{n+1}$ is a set of finite perimeter. Let $\Omega_0\subset\Omega$ be a  Lipschitz subdomain such that $\sigma(\partial^* \Omega \cap \partial \Omega_0)>0$. Let $\nu$ be the measure theoretic outer normal to $\Omega$, and $\nu_0$ be the measure theoretic outer normal to $\Omega_0$. Then  $\nu(\xi)=\nu_0(\xi)$ for  $\sigma$-almost every $\xi\in\partial^*\Omega\cap\partial\Omega_0$.
\end{lemma}

\subsection{The Haj\l{}asz-Sobolev space}\label{sec.sobolev}    Let $(\Sigma,\sigma)$ be a metric space with $\sigma$ a doubling
measure on $\Sigma$, which means that there is a uniform constant $C_\sigma\geq1$ such that $\sigma(B(x,2r))\leq C_\sigma\, \sigma(B(x,r))$, for all $x\in \Sigma$ and $ r>0$. 

For a Borel function $f:\Sigma\to\R$, we say that a non-negative Borel function $g:\Sigma \to \R$ is a {\it Haj\l asz upper gradient of  $f$} if  
\begin{equation}\label{eq:H-Sobolev}\nonumber
	|f(x)-f(y)| \leq |x-y| \,(g(x)+g(y))\quad \mbox{ for $\sigma$-a.e.\ $x, y \in \Sigma$.} 
\end{equation}
We denote the collection of all the Haj\l asz upper gradients of $f$ by $D(f)$. For $p\geq1$, we denote by $\dt{W}^{1,p}(\Sigma)$ the space of Borel functions $f$ which have 
a Haj\l asz upper gradient in $L^p(\sigma)$, and we let $W^{1,p}(\Sigma)$ be the space of functions $f\in L^p(\sigma)$ which have a Haj\l asz upper gradient in $L^p(\sigma)$; that is,  $W^{1,p}(\Sigma)=  \dt W^{1,p}(\Sigma) \cap L^p(\sigma)$.
We  define the semi-norm  
\begin{equation}\label{eqseminorm}\nonumber
	\| f \|_{\dt W^{1.p}(\Sigma)} = \inf_{g \in D(f)} \| g\|_{L^p(\Sigma)}
\end{equation}
and the scale-invariant norm
\begin{equation}\label{eqnorm}\nonumber
	\| f\|_{W^{1,p}(\Sigma)} = \diam(\Sigma)^{-1} \|f\|_{L^p(\Sigma)} +   \inf_{g \in D(f)} \| g\|_{L^p(\Sigma)}.
\end{equation}
For any   metric space $\Sigma$, if $p\in (1,\infty)$, from the uniform convexity of $L^p(\sigma)$, one easily deduces
that the infima in the definitions of   $\|\cdot\|_{W^{1,p}(\Sigma)}$ and $\|\cdot\|_{\dt W^{1,p}(\Sigma)}$   are uniquely attained. We denote by $\nabla_{H,p} f$ the function $g$ which attains the infimum. The spaces $\dt{W}^{1,p}(\Sigma)$ and ${W}^{1,p}(\Sigma)$ were introduced in \cite{haj96}, and are known as the {\it Haj\l{}asz-Sobolev spaces}.

\subsection{Non-tangential maximal functions and Carleson functions}\label{sec.ntmax}

For $\alpha>0$   and $\xi\in\partial\Omega$, we define the \emph{cone with vertex $\xi$ and aperture} $\alpha>0$ by
\begin{equation}\label{eq.cone}\nonumber
	\Gamma_\alpha^{\Omega}(\xi):=\{x\in\Omega:|x-\xi|<(1+\alpha)\delta(x)\}.
\end{equation}
If the domain or aperture are implicit from the surrounding context, then we will often omit them in the sequel for ease of notation. Define the \emph{non-tangential maximal function} of  $u\in L^{\infty}_{\loc}(\Omega)$ by 
\begin{equation}\label{eq.ntmax}
\m N_\alpha(u)(\xi):=\sup_{x\in\Gamma_\alpha(\xi)}|u(x)|,\quad\text{for }\xi\in\partial\Omega.
\end{equation}

Following \cite{kp}, we introduce the \emph{modified non-tangential maximal function} $\widetilde{\m N}_{\alpha,\hat c, r}$ for a given aperture $\alpha>0$, a parameter $\hat c\in(0,1/2]$, and $r\geq1$: for any $u\in L^r_{\loc}(\Omega)$, we write
\begin{equation}\label{eq.mnt}
	\widetilde{\m N}_{\alpha,\hat c,r}(u)(\xi):=\sup_{x\in\Gamma_\alpha(\xi)}\Big(\dashint_{B(x,\hat c\delta(x))}|u(y)|^r\,dm(y)\Big)^{1/r},\qquad\xi\in\partial\Omega.
\end{equation}
For $p\in(1,\infty)$, the $L^p$ norms of non-tangential maximal functions are equivalent under changes of   $\alpha$ or the averaging parameter $\hat c$. For ease of notation, we will often write   $\widetilde{\m N}_{\alpha,r}=\widetilde{\m N}_{\alpha,\hat c,r}$ if we do not wish to specify $\hat c$. When we do not need to specify neither $\alpha$ nor $\hat c$,  we will write ${\m N}=\m N_\alpha$,  and $\widetilde{\m N}_r=\widetilde{\m N}_{\alpha,\hat c,r}$. 

Let $q\in[1,\infty]$ and $\hat c\in(0,1/2]$.   Define the \emph{$q$-Carleson functional} of a   function $H:\Omega\ra\bb R, H\in L_{\loc}^{q}(\Omega)$ by
\begin{equation}\label{eq.carlesonfn}
	\n C_{\hat c,q}(H)(\xi):=\sup_{r>0}\frac1{r^n}\int_{B(\xi,r)\cap\Omega}\Big(\dashint_{B(x,\hat c\delta(x))}|H|^{q}\Big)^{1/q}\,dm(x),\qquad \xi\in\partial\Omega.
\end{equation} 
The case $q=\infty$ is defined with the appropriate analogue. For $p\in(1,\infty)$, the $L^p$ norms of the Carleson functionals $\n C_{\hat c,q}$ defined in (\ref{eq.carlesonfn}) are equivalent under a change of the averaging parameter $\hat c$, and thus we    write $\n C_q=\n C_{\hat c,q}$ if we do not need to specify $\hat c$. 

\begin{definition}[The spaces ${\bf N}_{r,p}$ and ${\bf C}_{r,p}$]\label{def.spaces} Fix a   domain $\Omega$ satisfying the corkscrew condition and with $n$-Ahlfors regular boundary. By $C_c(\overline{\Omega})$ we denote the space of continuous functions on $\overline\Omega$ (the closure of $\Omega$) with compact support in $\overline\Omega$. Fix $\alpha_0>0$, $\hat c_0\in(0,\frac12]$, and write $\wt{\m N}_r:=\wt{\m N}_{\alpha_0,\hat c_0,r}$, and $\n C_r:=\n C_{\hat c_0,r}$. Note that the functionals
\[
\Vert u\Vert_{{\bf N}_{r,p}}:=\Vert\wt{\m N}_r(u)\Vert_{L^p(\partial\Omega,\sigma)} ,\qquad\qquad\Vert u\Vert_{{\bf C}_{r,p}}:=\Vert\n C_r(u)\Vert_{L^p(\partial\Omega,\sigma)}
\]
define norms on $C_c(\overline{\Omega})$ and on $C_c(\Omega)$, respectively.  Let ${\bf N}_{r,p}$ be the completion of $C_c(\overline\Omega)$ under the norm $\Vert\cdot\Vert_{{\bf N}_{r,p}}$, and let ${\bf C}_{r,p}$ be the completion of $C_c(\Omega)$ under the norm   $\Vert\cdot\Vert_{{\bf C}_{r,p}}$. By \cite[Lemma 2.5]{mpt25}, if $r,p\in(1,\infty)$, we may also identify ${\bf C}_{r,p}$ as the Banach space of functions $u\in L^r_{\loc}(\Omega)$ such that $\Vert\n C_r(u)\Vert_{L^p(\partial\Omega)}<+\infty$. 
\end{definition}

The following proposition establishes a certain quantitative quasi-duality between the spaces ${\bf N}_{r,p}$ and ${\bf C}_{r',p'}$.   When $\Omega$ is the half-space, these estimates were shown by Hyt\"onen and Ros\'en in \cite[Theorem 3.1, Theorem 3.2]{hr13}. The general case cited here was shown in   \cite{mpt25}.

\begin{proposition}\label{prop.duality} Let $\Omega\subset\bb R^{n+1}$, $n\geq1$, be a domain satisfying the corkscrew condition and such that $\partial\Omega$ is $n$-Ahlfors regular. Suppose that either $\Omega$ is bounded, or that $\partial\Omega$ is unbounded. Let $p,q\in(1,\infty)$ and $p'$, $q'$ their H\"older conjugates. Then ${\bf N}_{q,p}=({\bf C}_{q',p'})^*$, and moreover,
	\begin{equation}\label{eq.hr7}\nonumber
		\Vert uH\Vert_{L^1(\Omega)}\lesssim\Vert\wt{\m N}_q(u)\Vert_{L^p(\partial\Omega)}\Vert\n C_{q'}H\Vert_{L^{p'}(\partial\Omega)},\qquad u\in L^q_{\loc}(\Omega), H\in L^{q'}_{\loc}(\Omega),
	\end{equation}
	\begin{equation}\label{eq.hr8}\nonumber
		\Vert\widetilde{\m N}_q(u)\Vert_{L^p(\partial\Omega)}\lesssim\sup_{H:\Vert\n C_{q'}(H)\Vert_{L^{p'}(\partial\Omega)}=1}\Big|\int_{\Omega}Hu\,dm\Big|,\qquad u\in L^q_{\loc}(\Omega),
	\end{equation}
	\begin{equation}\label{eq.hr9}\nonumber
		\Vert\n C_{q'}H\Vert_{L^{p'}(\partial\Omega)}\lesssim\sup_{u:\Vert\wt{\m N}_qu\Vert_{L^p(\partial\Omega)}=1}\Big|\int_\Omega Hu\,dm\Big|,\qquad H\in L^{q'}_{\loc}(\Omega).
	\end{equation}
\end{proposition}

The following extension result is fundamental to our program.

\begin{lemma}[The Varopoulos extension, {\cite[Theorem 1.4]{mz25}}]\label{lm.varopoulos} Let $\Omega\subset\bb R^{n+1}$, $n\geq2$, be a    domain satisfying the corkscrew condition and with $n$-Ahlfors regular boundary, and fix $g\in \Lip_c(\partial\Omega)$. Then there exists a function $\Phi:\overline{\Omega}\ra\bb R$ with the following properties:
	\begin{enumerate}[(i)]
		\item\label{item.v1} $\Phi\in C^{\infty}(\Omega)\cap\Lip(\overline{\Omega})\cap\dt W^{1,2}(\Omega)$.
		\item\label{item.v2} $\Vert\m N(\Phi)\Vert_{L^p(\sigma)}+\Vert\n C_\infty(\nabla\Phi)\Vert_{L^p(\sigma)}\lesssim\Vert g\Vert_{L^p(\sigma)}$, for each $p\in(1,\infty]$.
		\item\label{item.v3} $\Vert\m N(\delta_\Omega\nabla\Phi)\Vert_{L^p(\sigma)}\lesssim\Vert g\Vert_{L^p(\sigma)}$.
		\item\label{item.v4} $\Phi|_{\partial\Omega}=g$ continuously.
	\end{enumerate}
	We call $\Phi$ the \emph{Varopoulos extension} of $g$ onto $\Omega$.
\end{lemma}

\subsection{Elliptic PDE preliminaries}\label{sec.pde}
We assume throughout that $A$ is a real, not necessarily symmetric $(n+1)\times(n+1)$ matrix of   merely  bounded measurable coefficients in $\Omega$ verifying the strong ellipticity conditions 
\begin{equation}\label{eq.elliptic}
	\lambda|\xi|^2\leq\sum_{i,j=1}^{n+1}A_{ij}(x)\xi_i\xi_j,\qquad \Vert A\Vert_{L^{\infty}(\Omega)}\leq\frac1{\lambda},\quad x\in\Omega,\quad\xi\in\bb R^{n+1}.
\end{equation}  
Define the elliptic operator $L$ acting formally on real-valued functions $u$ by
\[
Lu=-\dv(A\nabla u)=-\sum_{i,j=1}^{n+1}\frac{\partial}{\partial x_i}\Big(a_{ij}\frac{\partial u}{\partial x_j}\Big).
\]
We write $A^T$ for the transpose of $A$, and $L^*=-\dv A^T\nabla$.

We will also make use of the conormal derivative. Let $\m T(\partial\Omega)$ denote the trace space of $W^{1,2}(\Omega)$; that is, $\m T(\partial\Omega)$ is the quotient space $W^{1,2}(\Omega)/W^{1,2}_0(\Omega)$.

\begin{definition}[Conormal Derivative]\label{def.conormal} Given $H\in\Lip_c(\Omega)$ and the unique weak solution $v\in W_0^{1,2}(\Omega)$ to the Poisson  problem
	\begin{equation}\label{eq.poissonh}
		\left\{\begin{aligned}-\dv A\nabla v&=H,\qquad&\text{in }&\Omega,\\ v&=0,\qquad&\text{on }&\partial\Omega.\end{aligned}\right. 
	\end{equation}
	define the functional
	\begin{equation}\label{eq.funct}
		\ell_v(\varphi)=B[v,\Phi]:=\int_\Omega A\nabla v\nabla\Phi\,dm-\int_\Omega H\Phi\,dm,\qquad\varphi\in\m T(\partial\Omega),
	\end{equation}
	where $\Phi\in W^{1,2}(\Omega)$ satisfies $\operatorname{Tr}\Phi=\varphi$.  It is easy to see that $\ell_v(\varphi)$ is well-defined and independent of the extension $\Phi$, since $v$ solves the equation $Lv=H$ in the weak sense. We call $\ell_v(\varphi)$ the \emph{conormal derivative of $v$}, and denote $\partial_{\nu_A}v=\ell_v$.
\end{definition}

\begin{lemma}[De Giorgi-Nash-Moser estimate]\label{lm.DGNM}
    Let \(B(x,2r)\subset\Omega\) and \(u\in W^{1,2}_{\mathrm{loc}}(\Omega)\) a weak solution to \(Lu=0\). Then
    \begin{align}\sup_{B(x,r)} |u|\lesssim \fint_{B(x,2r)} |u| \,dm.\label{eq.DGNM}\end{align}
\end{lemma}

\begin{lemma}[Boundary H\"{o}lder, {\cite[Corollary 4.49]{HS25}}]\label{lm.boundholder}
 Let \(\Omega\subset\mathbb{R}^{n+1}\) be a uniform domain with \(n\)-Ahlfors regular boundary. Let \(\xi\in\partial\Omega\), \(r\in (0,\mathrm{diam}(\Omega)/2)\). Suppose that \(u\in W^{1,2}_{\mathrm{loc}}(\Omega)\) is a weak solution of \(Lu = 0\) in \(\Omega\cap B(\xi,2r)\) with vanishing Neumann data on \(B(\xi, 2r)\cap\partial\Omega\).
 Then there exist \(\alpha>0\) and \(C>0\), depending only on \(n,\lambda\), the Ahlfors regularity constant, the Harnack chain constant and the corkscrew constant, such that for all \(x,y\in B(\xi,r)\)
    \[ |u(x)-u(y)|\lesssim \Big(\frac{|x-y|}{r}\Big)^\alpha \Big(\fint_{B(\xi,2r)} |u|^2 \,dm\Big)^{\frac{1}{2}}.\]

\end{lemma}

\begin{definition}[The $L^p$ Dirichlet problem $(\operatorname{D}_p^L)$]\label{def.dirichlet} Let $p\in(1,\infty)$. We say that the (homogeneous) \emph{Dirichlet problem for the operator $L$ with $L^{p}$ data in $\Omega$ is solvable} (write $(\Di_{p}^L)$ is solvable in $\Omega$), if there exists $C\geq1$ so that for each $g\in C(\partial\Omega)$, the solution $u$ to the continuous Dirichlet problem 
\begin{equation}\label{eq.dirichlet2} 
\left\{\begin{aligned}-\dv A\nabla u&=0,\quad&\text{in }&\Omega,\\u&=g,\quad&\text{on }&\partial\Omega,\end{aligned}\right.
\end{equation}
satisfies the estimate 
\begin{equation}\label{eq.direst}\nonumber
\Vert\m N(u)\Vert_{L^{p}(\partial\Omega,\sigma)}\leq C\,\Vert g\Vert_{L^{p}(\partial\Omega,\sigma)}.
\end{equation}
\end{definition}

\begin{definition}[The $L^p$ Regularity problem $(\operatorname{R}_p^L)$]\label{def.regularity} Let   $p\in(1,\infty)$. We say that the (homogeneous) \emph{Dirichlet regularity problem}   (or just \textit{regularity problem})  for the operator $L$ with   $\dt W^{1,p}(\pom,\sigma)$  data  is solvable in $\Omega$ (write $(\Reg_{p}^L)$ is solvable in $\Omega$), if there exists $C\geq1$ so that for each $f\in\Lip(\partial\Omega)$, the solution $u$ to the continuous Dirichlet problem (\ref{eq.dirichlet2}) with boundary data $f$ satisfies the estimate
\begin{equation}\label{eq.regest}
	\Vert\wt{\m N}_2(\nabla u)\Vert_{L^{p}(\partial\Omega,\sigma)}\leq C\Vert f\Vert_{\dt{W}^{1,p}(\partial\Omega,\sigma)}.
\end{equation}
\end{definition}


Let $\Omega\subset\bb R^{n+1}$ be a bounded chord-arc domain.  We say that the \emph{variational Neumann problem} is solvable in $\Omega$ if for each $g\in L^{\frac{2n}{n+1}}(\partial\Omega)$ with $\int_{\partial\Omega}g\,d\sigma=0$, there exists a solution $u\in W^{1,2}(\Omega)$ to the problem
\begin{equation}\label{eq.varneumann} 
\left\{\begin{aligned}-\dv A\nabla u&=0,\quad&\text{in }&\Omega,\\ \partial_{\nu_A}u&=g,\quad&\text{on }&\partial\Omega.\end{aligned}\right.
\end{equation}
Note that the variational solution is unique up to a constant. For more details on the solvability on the variational Neumann problem, please see \cite[Section 2.8]{mt26}.

\begin{definition}[The $L^p$ Neumann problem $(\operatorname{N}_p^L)$]\label{def.neumann} Let $\Omega\subset\bb R^{n+1}$ be a bounded uniform domain with $n$-Ahlfors regular boundary and fix $p\in(1,\infty)$. We say that the \emph{Neumann problem for the operator $L$ with $L^p$ data in $\Omega$ is solvable} (write $(\operatorname{N}_{p}^L)$ is solvable in $\Omega$), if there exists $C\geq1$ so that for each $g\in L^p(\partial\Omega)\cap L^{\frac{2n}{n+1}}(\partial\Omega)$ with $\int_{\partial\Omega}g\,d\sigma=0$, the variational solution to (\ref{eq.varneumann})   satisfies the estimate
\begin{equation}\label{eq.neumannest}
\Vert\wt{\m N}_2(\nabla u)\Vert_{L^{p}(\partial\Omega,\sigma)}\leq C\Vert g\Vert_{L^p(\partial\Omega,\sigma)}.
\end{equation}

\end{definition}


Recall that $\widehat{W}^{1,2}(\Omega)$ denotes the space of functions in $W^{1,2}(\Omega)$ with zero mean value of the trace on $\partial\Omega$.
\begin{definition}[The rough $L^p$ Neumann problem]\label{def.roughneumann} Let $\Omega\subset\bb R^{n+1}$ be a bounded uniform domain with $n$-Ahlfors regular boundary and fix $p\in(1,\infty)$. We say that the \emph{rough Neumann problem for the operator $L$ with $L^p$ data in $\Omega$ is solvable}, if there exists $C\geq1$ so that for all $j,k\in\{1,\ldots,n+1\}$ and each $g\in \Lip(\partial_*\Omega)$, the variational solution $u\in\widehat W^{1,2}(\Omega)$ of 
\begin{equation}\label{eq.roughneumann} 
\left\{\begin{aligned}-\dv A\nabla u&=0,\quad&\text{in }&\Omega,\\ \partial_{\nu_A}u&=\partial_{t,j,k}g,\quad&\text{on }&\partial^*\Omega.\end{aligned}\right.
\end{equation}
satisfies the estimate
 \begin{equation}\label{eq.roughneumannest}
 \Vert\m N(u)\Vert_{L^{p}(\partial\Omega,\sigma)}\leq C\Vert g\Vert_{L^p(\partial^*\Omega,\sigma)}.
 \end{equation}

 \end{definition}

\section{Failure of Approximation Theorem for the Dirichlet problem at high orders of regularity}\label{sec.failure}

In this  section we give the short proofs of the propositions concerning the failure of the approximation theorem at the levels of regularity of the variational Dirichlet problem and the Dirichlet-Regularity problem.

    
\noindent\emph{Proof of Proposition \ref{prop.approxReg}.} Let $\tilde p:=\min\{p,2\}$. Since $\Omega$ is bounded, then there exists a constant $C_\Omega$ such that
\begin{equation}\label{eq.bound}
\Vert\nabla w_j\Vert_{L^{\tilde p}(\Omega)}\leq C_\Omega\Vert\m N_2(\nabla w_j)\Vert_{L^p(\partial\Omega)},\qquad\text{for all }j,
\end{equation}
where $C_\Omega$ depends on $\diam\Omega$. Moreover, since $\{\nabla w_j\}_j$ is a weak$-*$ convergent sequence, it is in particular uniformly bounded in ${\bf N}_{2,p}$ by the Uniform Boundedness Principle.
Consequently, the sequence $\{\nabla w_j\}_j$ is uniformly bounded in $L^{\tilde p}(\Omega)$ by (\ref{eq.bound}). Since each of the $w_j$ has a vanishing trace (in the sense of belonging to $W_0^{1,2}(\Omega)$), then by the Poincar\'{e} inequality we conclude that $\{w_j\}_j$ is uniformly bounded in $W_0^{1,\tilde p}(\Omega)$.  

By the Rellich-Kondrachov Theorem, it follows that there exists a subsequence $\{w_{j_k}\}_k$ that converges weakly in $W^{1,\tilde p}(\Omega)$ and strongly in $L^{\tilde p}(\Omega)$. But since $W_0^{1,\tilde p}(\Omega)$ is a closed subspace of $W^{1,\tilde p}(\Omega)$, $w_{j_k}$ converges weakly in $W_0^{1,\tilde p}(\Omega)$ and strongly in $L^{\tilde p}(\Omega)$. Let us denote this limit function by $w\in W_0^{1,\tilde p}(\Omega)$.

    Since $\Lip_c(\Omega)\subset L^{\tilde p'}(\Omega)\cap{\bf C}_{2,p'}(\Omega)$, we obtain for every $F\in (\Lip_c(\Omega))^{n+1}$ that
    \begin{align*}
        \int_\Omega (\nabla u-\nabla w) F \, dm= \int_\Omega (\nabla u-\nabla w_{j_k}) F \, dm + \int_\Omega (\nabla w_{j_k}-\nabla w) F \, dm\lra 0
    \end{align*}
    as $k\to \infty$. Due to density from $\Lip_c(\Omega)$ in $L^{\tilde p'}(\Omega)$, we conclude that $\Vert \nabla u-\nabla w\Vert_{L^{\tilde p}}=0$, and hence $u=w+c$ for some constant $c$. But since $w$ has vanishing trace, we get that $f\equiv u|_{\partial\Omega}\equiv w|_{\partial\Omega}+c\equiv c$, as desired.\hfill{$\square$}

\noindent\emph{Proof of Proposition \ref{prop.approxvar}.} This proof proceeds very similarly to that of Proposition \ref{prop.approxReg}. Since $\{\nabla w_j\}_j$ is a weakly convergent sequence in $L^2(\Omega)$, then it is in particular a bounded sequence in $L^2(\Omega)$.  Since each of the $w_j$ has a vanishing trace (in the sense of belonging to $W_0^{1,2}(\Omega)$), then by the Poincar\'{e} inequality we conclude that $\{w_j\}_j$ is uniformly bounded in $W_0^{1,2}(\Omega)$.  

It follows that there exists a subsequence $\{w_{j_k}\}_k$ that converges weakly in $W^{1,2}(\Omega)$. But since $W_0^{1,2}(\Omega)$ is a closed subspace of $W^{1,2}(\Omega)$, $w_{j_k}$ converges weakly in $W_0^{1,2}(\Omega)$. Let us denote this limit function by $w\in W_0^{1,2}(\Omega)$. Since $\nabla u$ is the weak limit of $\{\nabla w_j\}_j$ in $L^2(\Omega)$, we must have that $\nabla u\equiv\nabla w$ pointwise almost everywhere in $\Omega$, and thus $u=w+c$   for some constant $c$. But since $w$ has vanishing trace, we get that $f\equiv u|_{\partial\Omega}\equiv w|_{\partial\Omega}+c\equiv c$, as desired.\hfill{$\square$}

\section{Approximation for the Neumann problem}\label{sec.neumann}

In this section, we show Theorem \ref{thm.approxNeumann}.

\noindent\emph{Proof of Theorem \ref{thm.approxNeumann}.}

\emph{Proof of \ref{item.c1}.} Suppose that $g\in\Lip(\partial\Omega)$ with $\int_{\partial\Omega}g\,d\sigma=0$, and let $\Phi=\Phi_g$ be the Varopoulos  extension of $g$ onto $\Omega$ as in Lemma \ref{lm.varopoulos}. Now, for each $j\in\bb N$ large enough (namely, such that $2^{-j}<\diam (\Omega)/2$), let $\eta_j$ be a cut-off function satisfying that $\eta_j\equiv1$ on $B(\partial\Omega,2^{-(j+1)})$, $\supp\eta_j\subset B(\partial\Omega,2^{-j})$, $|\nabla\eta_j|\lesssim2^{j}$, and $0\leq\eta_j\leq1$. We set 
\begin{equation}\label{eq.fj0Neu}
\tilde H_j:=\rho_j\Phi\eta_j-\fint_\Omega \rho_j\Phi\eta_j \,dm,
\end{equation}
where $\rho$ is given by
\begin{equation}\label{eq.rhoNeu}
\rho_j(x):=\int_{\partial\Omega}\Big(\int_{\Gamma_\alpha(\xi)}\eta_j\,dm\Big)^{-1}\chi_{\Gamma_\alpha(\xi)}(x)d\sigma(\xi).
\end{equation}
By the corkscrew condition, we can see that $\int_{\Gamma_\alpha(\xi)}\eta_j\,dm\approx 2^{-j(n+1)}$, and thus, Ahlfors regularity of the boundary gives $\rho_j(x)\approx 2^{j(n+1)}\delta(x)^{n}$ with uniform constants.We claim that
\begin{align}
\Vert\n C_\infty(\rho_j\Phi\eta_j)\Vert_{L^p(\sigma)}\lesssim\Vert \m M(\m N(\Phi))\Vert_{L^p(\sigma)}\lesssim\Vert g\Vert_{L^p(\sigma)},\label{eq:claimCBoundedByNNeu}
\end{align}
where the implicit constant depends on the aperture \(\alpha\), on the corkscrew constant, and on the Ahlfors regularity constant. We postpone its proof until the end of part 1. 

Hence $\tilde H_j\in{\bf C}_{\infty,p}$, and by H\"older's inequality in fact we have that $\tilde H_j\in{\bf C}_{q,p}$ for all $q\in(1,\infty]$. Moreover, $\tilde{H}_j$ has zero mean-value. Now, for each $q\in(1,\infty)$, by the density of $\Lip_c(\Omega)$ in ${\bf C}_{q,p}$ (see Lemma~\ref{lem:patch-density}), we may find compactly supported Lipschitz functions $H_j$ with zero mean-value such that
\begin{equation}\label{eq.fjNeu}
\Vert\n C_q(H_j-\tilde H_j)\Vert_{L^p(\sigma)}\leq2^{-j}\Vert g\Vert_{L^p(\sigma)}\qquad\textrm{and}\qquad \supp(H_j)\subset B(\partial\Omega, 2^{-j+1}).
\end{equation}

This finishes the proof of \ref{item.c1} when $g$ is Lipschitz continuous.

\emph{Proof of \ref{item.c2}.} We continue to assume that  $g\in\Lip(\partial\Omega)$. Since $\Lip_c(\Omega)\subset L^2(\Omega)$, and $H_j$ has zero mean value, Lax-Milgram yields solutions $w_j\in \widehat{W}^{1,2}(\Omega)$ to \eqref{eq.poissonNeumann0} with right hand side $H_j$.

Note also that the solvability of $(\N_p^L)$ implies that $\nabla u\in({\bf N}_{2,p})^{n+1}$. Let $p'$ be the H\"older conjugate of $p$. By Proposition~\ref{prop.duality}, we know that $({\bf C}_{2,p'})^*={\bf N}_{2,p}$, and thus $(({\bf C}_{2,p'})^{n+1})^*=({\bf N}_{2,p})^{n+1}$. Hence, it suffices to show

\begin{equation}\label{eq.convNeuNeumann}
\int_\Omega\nabla(w_j-u)F\,dm\lra0,\quad\text{as }j\ra\infty, \text{ for each }F\in({\bf C}_{2,p'})^{n+1}.
\end{equation} 

Assume first that $F\in(\Lip_c(\Omega))^{n+1}$, and let $V\in \widehat{W}^{1,2}(\Omega)$ be the weak solution to the Poisson problem
\begin{equation}\label{eq.poissonregularityNeumann}
	\left\{\begin{aligned}-\dv A^T\nabla V&=-\dv F,\qquad&\text{in }&\Omega,\\ \partial_{\nu_{A^T}}V&=0,\qquad&\text{on }&\partial\Omega.\end{aligned}\right. 
\end{equation}
It follows that
\begin{equation}\label{eq.conv2}
	\int_\Omega \nabla (w_j-u)F\,dm=\int_\Omega A^T\nabla V\nabla(w_j-u)\,dm,
\end{equation}
where we used that $\partial_{\nu_{A^T}}V\equiv0$.
On the other hand, note that
\begin{multline}\label{eq.conv3}
	\int_\Omega A^T\nabla V\nabla(w_j-u)\,dm=\int_\Omega A^T\nabla V\nabla w_j\,dm-\int_\Omega A^T\nabla V\nabla u\,dm\\ =\int_\Omega A\nabla w_j\nabla V\,dm-\int_\Omega A\nabla u\nabla V\,dm=\int_\Omega H_j V\,dm-(\partial_{\nu_A}u,V),
\end{multline}
since $-\mathrm{div} A\nabla u=0$, $\partial_{\nu_A}w_j\equiv 0$, $-\dv A\nabla w_j= H_j$ in $\Omega$ and $V\in \widehat{W}^{1,2}(\Omega)$. Next, using \eqref{eq.fj0Neu}, \eqref{eq.fjNeu}, and that the mean value of $V$ is zero, we have
\begin{equation}\label{eq.conv31}
\int_\Omega H_j V\,dm=\int_\Omega(H_j-\tilde H_j) V\,dm+\int_\Omega \rho_j\Phi\eta_j V\,dm.
\end{equation}

Putting (\ref{eq.conv2}) together with  (\ref{eq.conv3}), (\ref{eq.conv31}), we get that
\begin{equation}\label{eq.conv4}
	\int_\Omega \nabla(w_j-u)F\,dm=\int_\Omega(H_j-\tilde H_j) V\,dm+\int_\Omega \rho_j\Phi\eta_j V\,dm~-~(\partial_{\nu_{A}}u,V).
\end{equation}
By Proposition~\ref{prop.duality} we obtain for the first term
\begin{equation}\label{eq.fjapprox}
\Big|\int_\Omega(H_j-\tilde H_j)V\,dm\Big|\lesssim\Vert\n C_2(H_j-\tilde H_j)\Vert_{L^p(\sigma)}\Vert\widetilde{\m N}_2(V)\Vert_{L^{p'}(\sigma)}\lesssim2^{-j}\Vert g\Vert_{L^p(\sigma)}\Vert\n C_2(F)\Vert_{L^{p'}(\sigma)},
\end{equation}
where we used that the Poisson-Neumann problem $(\operatorname{PN}_{p'}^{L^*})$ is solvable in $\Omega$, since $(\operatorname{D}_p^L)$ and $(\operatorname{N}_p^L)$ are solvable in $\Omega$ (see \cite[Definition 1.11, Lemma 7.1, and Lemma 3.1]{fl24}).

The second and third term of \eqref{eq.conv4} can be combined. First, we use the definition of $\rho_j$ (see \eqref{eq.rhoNeu}) and Fubini-Tonelli to rewrite the second term as
\begin{align*}
    \int_\Omega \rho\Phi\eta_j V\,dm
    =\int_{\Omega} \left[\int_{\partial\Omega}\Big(\int_{\Gamma_\alpha(\xi)}\eta_j\,dm\Big)^{-1}\chi_{\Gamma_\alpha(\xi)}(x)\,d\sigma(\xi)\right]\Phi(x)\eta_j(x) V(x)\,dm(x)
    \\
    =\int_{\partial\Omega}\left[\int_\Omega\Big(\int_{\Gamma_\alpha(\xi)}\eta_j\,dm\Big)^{-1}\chi_{\Gamma_\alpha(\xi)}(x)\Phi(x)\eta_j(x) V(x)\,dm(x)\right]\,d\sigma(\xi)
    \\
    =\int_{\partial\Omega}\left[\int_{\Gamma_\alpha(\xi)}\Big(\int_{\Gamma_\alpha(\xi)}\eta_j\,dm\Big)^{-1}\Phi(x)\eta_j(x) V(x)\,dm(x)\right]\,d\sigma(\xi).
\end{align*}
We abbreviate notation by setting \(m_j(x):=\Big(\int_{\Gamma_\alpha(\xi)}\eta_j\,dm\Big)^{-1}\eta_j(x)\) and noticing that \(\supp(m_j)\subset B(\partial\Omega,2^{-j})\) with \(\int_{\Gamma_\alpha(\xi)}m_j(x)dx=1\). Now, combining the second and third term of \eqref{eq.conv4} gives

\begin{align}
    \int_{\partial\Omega}\left[\int_{\Gamma_\alpha(\xi)}m_j(x)\Phi(x) V(x)\,dm(x)\right]\,d\sigma(\xi)-\int_{\partial\Omega} g(\xi)V(\xi)\,d\sigma(\xi)\nonumber
    \\
    =\int_{\partial\Omega}\left[\int_{\Gamma_\alpha(\xi)}m_j(x)\big(\Phi(x) V(x)-g(\xi)V(\xi)\big)\,dm(x)\right]\,d\sigma(\xi).\label{eq.conv5}
\end{align}
For $x\in \mathrm{supp}(m_j)\cap\Gamma_\alpha(\xi)$, we use triangle inequality to bound 
\begin{align}|\Phi(x) V(x)-g(\xi)V(\xi)|\leq |g(\xi)||V(x)-V(\xi)|+|\Phi(x)-g(\xi)||V(x)|.\label{eq.1}
\end{align}
By Lemma~\ref{lm.varopoulos}, \(\Phi\) is Lipschitz as the Varopolous extension of $g$, and we have 
\begin{align}
|\Phi(x)-g(\xi)|\leq L2^{-j},\label{eq.2}
\end{align}
where $L$ is the Lipschitz constant of $\Phi$. Furthermore, there exists a sufficiently large \(j_0\) such that $B(\partial\Omega, 2^{-j+6})\cap \mathrm{supp}(F)= \emptyset$ for all \(j\geq j_0\), where \(F\) is the compactly supported Poisson data in the PDE for \(v\) (see \eqref{eq.poissonregularityNeumann}). Thus, if \(j\geq j_0\) and \(x\in \Gamma_\alpha(\xi)\cap\mathrm{supp}(m_j)\), then Lemma~\ref{lm.boundholder} implies \begin{align}
|V(\xi)-V(x)|\lesssim  2^{-\beta (j-j_0)}\frac{1}{m(B(\xi,2^{-j_0})\cap\Omega)^{1/2}}\Big(\int_{B(\xi,2^{-j_0})\cap\Omega}|V|^2 \,dm\Big)^{1/2}\leq C 2^{-\beta(j-j_0)},\label{eq.3}
\end{align}
where the constant \(C\) depends on \(j_0\), the corkscrew constant of \(\Omega\) and \(\Vert V\Vert_{W^{1,2}(\Omega)}\).

Hence, combining \eqref{eq.1}, \eqref{eq.2}, and \eqref{eq.3}, we obtain for \eqref{eq.conv5}
\begin{align*}
    &\left|\int_{\partial\Omega}\left[\int_{\Gamma_\alpha(\xi)}m_j(x)\big(\Phi(x) V(x)-g(\xi)V(\xi)\big)\,dm(x)\right]\,d\sigma(\xi)\right|
    \\
    &\lesssim \int_{\partial\Omega}\left[\int_{\Gamma_\alpha(\xi)}m_j(x)\big(|g(\xi)||V(x)-V(\xi)|+|\Phi(x)-g(\xi)||V(x)|\big)\,dm(x)\right]\,d\sigma(\xi)
    \\
    &\lesssim \int_{\partial\Omega}\left[\int_{\Gamma_\alpha(\xi)}m_j(x)\big(2^{-\beta j}|g(\xi)|+2^{-j}|V(x)|)|\big)\,dm(x)\right]\,d\sigma(\xi)
    \\
    &\lesssim 2^{-\beta j}\Vert g\Vert_{L^1(\sigma)} + 2^{-j}\Vert \mathcal{N}^{2^{-j+2}}(V)\Vert_{L^1(\sigma)}.
\end{align*}
In the last line $\mathcal{N}^{2^{-j}}(V)(\xi):=\sup_{y\in\Gamma_\alpha(\xi)\cap B(\xi, 2^{-j+2})}|V(y)|$ is the truncated nontangential maximal function. Again, for \(j\geq j_0\), since the Poisson data of $V$ is supported away from \(B(\xi, 2^{-j+2})\), we can use Lemma~\ref{lm.DGNM} for $V$ on $B(\partial\Omega, 2^{-j+3})$, and hence we have 


\begin{align*}
    \mathcal{N}^{2^{-j+2}}(V)(\xi)&\leq\sup_{y\in\Gamma_\alpha(\xi)\cap B(\xi, 2^{-j+2})}|V(y)|
    \overset{\eqref{eq.DGNM}}{\lesssim}\sup_{y\in\Gamma_\alpha(\xi)\cap B(\xi, 2^{-j+2})}\fint_{B(y,\delta(y)/2)}|V(x)|\,dm(x) \leq \widetilde{\m N}_2(V)(\xi).
\end{align*} 
As in \eqref{eq.fjapprox}, by solvability of the Poisson-Neumann problem $(\operatorname{PN}_{p'}^{L^*})$, it follows that
\begin{align*}
     2^{-\beta j}\Vert g\Vert_{L^1(\sigma)} + 2^{-j}\Vert \mathcal{N}^{2^{-j+2}}(V)\Vert_{L^1(\sigma)}\lesssim  2^{-\beta j}\Vert g\Vert_{L^p(\sigma)} + 2^{-j}\Vert \widetilde{\mathcal{N}}_2(V)\Vert_{L^{p'}(\sigma)}\lesssim  2^{-\beta j}\Vert g\Vert_{L^p(\sigma)} + 2^{-j}\Vert \n C_2(F)\Vert_{L^{p'}(\sigma)},
\end{align*}
where the implicit constants also depend on the finite length of $\partial\Omega$.
In total with \eqref{eq.fjapprox}, we obtain for \(j\geq j_0\)
\[\left|\int_\Omega\nabla(w_j-u)F\,dm\right|\lesssim 2^{-j}\Vert g\Vert_{L^p(\sigma)}\Vert\n C_2(F)\Vert_{L^{p'}(\sigma)} + 2^{-\beta j}\Vert g\Vert_{L^p(\sigma)} + 2^{-j}\Vert \n C_2(F)\Vert_{L^{p'}(\sigma)}.\]
Taking the limit for \(j\to\infty\) concludes the proof of \eqref{eq.convNeuNeumann} for $H\in\Lip_c(\Omega)$.

Now suppose that $F\in({\bf C}_{2,p'})^{n+1}$ is arbitrary. Since compactly supported Lipschitz functions are dense in ${\bf C}_{2,p'}$ (see Lemma~\ref{lem:patch-density}), there is a sequence $\{F_k\}_{k\in\bb N}\subset(\Lip_c(\Omega))^{n+1}$ such that 
\[
\Vert\n C_2(F-F_k)\Vert_{L^{p'}(\sigma)}\ra0\qquad\text{as }k\ra\infty.
\]
We may write for each $j,k\in\bb N$,
\begin{equation}\label{eq.dense}
\int_\Omega \nabla(w_j-u)F\,dm=\int_\Omega \nabla(w_j-u)(F-F_k)\,dm~+~\int_\Omega\nabla(w_j-u)F_k\,dm=:I_1+I_2.
\end{equation}
Note that
\begin{align*}
|I_1|&\lesssim\Vert\widetilde{\m N}_2(\nabla (w_j-u))\Vert_{L^p(\sigma)}\Vert\n C_2(F-F_k)\Vert_{L^{p'}(\sigma)}\\ &\leq\Big[\Vert\widetilde{\m N}_2(\nabla w_j)\Vert_{L^p(\sigma)}+\Vert\widetilde{\m N}_2(\nabla u)\Vert_{L^p(\sigma)}\Big]\Vert\n C_2(F-F_k)\Vert_{L^{p'}(\sigma)}\\ &\lesssim\Big[\Vert\n C_2(H_j)\Vert_{L^p(\sigma)}+\Vert g\Vert_{L^p(\sigma)}\Big]\Vert\n C_2(F-F_k)\Vert_{L^{p'}(\sigma)}\lesssim\Vert g\Vert_{L^p(\sigma)}\Vert\n C_2(F-F_k)\Vert_{L^{p'}(\sigma)},
\end{align*}
where we used Proposition~\ref{prop.duality}, the solvability of $(\N_p^L)$, and the solvability of the Poisson-Neumann-Regularity problem $(\operatorname{PNR}_p^L)$, which follows from the solvability of $(\N_p^L)$ and $(\Di_{p'}^{L^*})$ (see \cite[Lemma 7.1]{fl24}). If we fix $\ep>0$, then the previous calculations show that we can fix $k$ large enough so that $|I_1|<\ep/2$ for all $j\in\bb N$. With $k$ fixed, since $H_k\in\Lip_c(\Omega)$, we have that for all $j$ large enough, $|I_2|<\ep/2$, as we already proved (\ref{eq.convNeuNeumann}) when $H\in\Lip_c(\Omega)$. We conclude that $\nabla w_j\ra \nabla u$ as $j\ra \infty$ in the weak$-*$ topology in $({\bf N}_{2,p})^{n+1}$. 

To complete the proof of \ref{item.c2} for $g\in \Lip(\partial\Omega)$, it remains to prove the claim \eqref{eq:claimCBoundedByNNeu}. Recall that \(j\) is large enough such that \(2^{-j}<\mathrm{diam}(\Omega)/2\). Fix $\zeta\in \partial\Omega$, then by \eqref{eq.rhoNeu} we have 
\begin{align*}
\n C_\infty(\rho_j\Phi\eta_j)(\zeta)&\leq \sup_{r>0}\frac{1}{r^n}\int_{B(\zeta, r)\cap B(\partial\Omega,2^{-j})}\sup_{y\in B(x,\delta(x)/2)}|\rho_j(y)\Phi(y)| dx
\\
&\lesssim \sup_{r>0}\frac{2^{-j(n+1)}}{r^n}\int_{B(\zeta, r)\cap B(\partial\Omega,2^{-j})}\delta(x)^n\sup_{y\in B(x,\delta(x)/2)}|\Phi(y)| dx,
\end{align*}
where the implicit constants depend on the aperture \(\alpha\), on the corkscrew constant, and on the Ahlfors regularity constant.
We observe that \cite[Lemma 3.13 and Lemma 3.14]{mpt} also hold on domains that merely satisfy the corkscrew condition and have Ahlfors regular boundary, and hence 
\begin{align*}
&\sup_{r>0}\frac{2^{j(n+1)}}{r^n}\int_{B(\zeta, r)\cap B(\partial\Omega,2^{-j})}\delta(x)^n\sup_{y\in B(x,\delta(x)/2)}|\Phi(y)| dx
\\
&\lesssim \sup_{r>0}\frac{2^{j(n+1)}}{r^n}\int_{B(\zeta, r)}\int_{\Gamma_\alpha(P)\cap B(\partial\Omega,2^{-j})}\sup_{y\in B(x,\delta(x)/2)}|\Phi(y)| \,dx \,d\sigma(P)
\\
&\lesssim \sup_{r>0}\frac{1}{r^n}\int_{B(\zeta, r)}\m N(\Phi) \,d\sigma=\m M(\m N(\Phi))(\zeta).
\end{align*}
Here we used in the last line that $m(\Gamma_\alpha(P)\cap B(\partial\Omega,2^{-j}))\lesssim 2^{-j(n+1)}$ for every $P\in\partial\Omega$. This completes the proof of the claim \eqref{eq:claimCBoundedByNNeu}.

{\bf Part 2. Case $g\in L^p(\partial\Omega,\sigma)$.}

Now, we fix $g\in L^p(\partial\Omega,\sigma)$ with $\int_{\partial\Omega} gd\sigma=0$, and let $k\in\mathbb{N}$. By \cite[Lemma 2.3]{mt26}, there exists $g_k\in\Lip(\partial\Omega)$ with $\int_{\partial\Omega} g_kd\sigma=0$ and $\Vert g_k-g\Vert_{L^p}<\frac{1}{k}$. Let us denote by $u_k$ the solution to \eqref{eq.neumann0} with boundary data $g_k$. Let \(\{F_l\}_{l\in\mathbb{N}}\subset(\Lip_c(\Omega))^{n+1}\) be a countable subset that is dense in \(({\bf C}_{2,p'})^{n+1}\) (see Lemma~\ref{lem:patch-separability}). Furthermore by Part 1, there exists a sequence $H_j^k\in\Lip_c(\partial\Omega)$ and a sequence of solutions $w_j^k\in \widehat{W}^{1,2}(\Omega)$ to \eqref{eq.poissonNeumann0} such that $\Vert H_j^k\Vert_{{\bf C}_{2,p}}\leq C\Vert g_k\Vert_{L^p}$ and $\nabla w_j^k\ra \nabla u_k$ as $j\ra\infty$ in the weak-$*$ topology in $({\bf N}_{2,p})^{n+1}$. For every \(k\), we can now choose \(j_k\) such that 
\begin{equation}\label{eq.flarg}
\Big|\int_\Omega (\nabla w_{j_k}^{k}-\nabla u_k)F_l \,dm\Big|< \frac{1}{k} \qquad \textrm{  for all }1\leq l\leq k.
\end{equation}

If we set \(H_k:=H_{j_k}^{k}\) and \(w_k:=w_{j_k}^{k}\), then we claim that the sequences $H_k$ and $w_k$ satisfy the assertions of the theorem. Indeed,
\[\Vert H_k\Vert_{{\bf C}_{2,p}}=\Vert H^{k}_{j_k}\Vert_{{\bf C}_{2,p}}\leq C\Vert g_k\Vert_{L^p}\leq C(\Vert g-g_k\Vert_{L^p}+\Vert g_k\Vert_{L^p})\leq 2C\Vert g\Vert_{L^p},\]
which shows \eqref{eq.uniformcNeu}. Moreover, for every \(\varepsilon>0\) and for every \(F\in ({\bf C}_{2,p'})^{n+1}\) there exists \(k_0=k_0(\varepsilon,F)\) such that \(\Vert F_{k_0}-F\Vert_{{\bf C}_{2,p'}}<\varepsilon\). Then,  for all \(k\geq k_0\), we have
\begin{align*}
    \int_\Omega (\nabla w_k-\nabla u)F\, dm&=\int_\Omega (\nabla w_k-\nabla u)(F-F_{k_0})\, dm + \int_\Omega (\nabla w_k-\nabla u)F_{k_0}\, dm=I_1+I_2.
\end{align*}
By \eqref{eq.flarg}, Proposition~\ref{prop.duality}, and the solvability of $(\operatorname{N}_p^L)$, we have for all $k\geq k_0$ that
\begin{multline*}|I_2|\leq \Big|\int_\Omega (\nabla w_{j_k}^{k}-\nabla u_{k})F_{k_0}\,dm\Big| + \Big|\int_\Omega (\nabla u_{k}-\nabla u)F_{k_0}\,dm\Big|\\ <\frac{1}{k}+\Vert \nabla(u_{k}-u)\Vert_{{\bf N}_{2,p}}\Vert F_{k_0}\Vert_{{\bf C}_{2,p'}}\leq \frac{1}{k}+C\Vert g_{k}-g\Vert_{L^p}\Vert F\Vert_{{\bf C}_{2,p'}}\leq\frac{C\Vert F\Vert_{{\bf C}_{2,p'}}}{k}.
\end{multline*}
By Proposition~\ref{prop.duality} and the solvability of \((\N_p^L)\) and of the Poisson-Neumann-Regularity problem $(\operatorname{PNR}_p^L)$, which follows from the solvability of $(\N_p^L)$ and $(\Di_{p'}^{L^*})$ (see \cite[Lemma 7.1]{fl24}), we get
\begin{align*}
    |I_1|&\leq(\Vert \nabla w_k\Vert_{{\bf N}_{2,p}}+\Vert \nabla u\Vert_{{\bf N}_{2,p}})\Vert F_{k_0}-F\Vert_{{\bf C}_{2,p'}}\lesssim \varepsilon(\Vert H_k\Vert_{{\bf C}_{2,p}}+\Vert g\Vert_{L^p})\lesssim \varepsilon\Vert g\Vert_{L^p}.
\end{align*}
This completes the proof of \(\nabla w_k\to \nabla u\) in the weak-* topology in \(({\bf N}_{2,p})^{n+1}\).

It remains to prove strong local convergence of the gradients.
Fix a ball $B$ with $\overline{2B}\subset\Omega$, and set
$v_k=w_k-u$. Since the $L^p$ norm of the nontangential maximal function of \(\nabla v_k\) is uniformly bounded in
$L^2(2B)$, the weak-$*$ convergence in \(({\bf N}_{2,p})^{n+1}\)
implies
\[
 \nabla v_k\rightharpoonup0\quad\text{in }L^2(2B).
\]
Indeed, for every \(F\in (L^2(2B))^{n+1}\) that is extended by zero out side of \(2B\) and hence lies in \({\bf C}_{2,p'}(\Omega)\) we have
\[\int_{2B}\nabla v_k\cdot F \,dm=\int_{\Omega}\nabla v_k\cdot F{\bf 1}_{2B} \,dm\to 0 \qquad \textrm{ as }k\to\infty.\]

Let us denote the mean value by \((v_k)_{2B}:=\int_{2B}v_k \, dm\).
By Poincar\'e's inequality, $\{v_k-(v_k)_{2B}\}$ is uniformly bounded in $W^{1,2}(2B)$.
Every subsequence has a further subsequence converging strongly in
$L^2(2B)$ by the Rellich--Kondrachov theorem. Every such limit has gradient zero, mean-value zero, and therefore equals zero. Hence
$v_k-(v_k)_{2B}\to0$ strongly in $L^2(2B)$ along the whole sequence. By construction of \(H_k\), $Lv_k=0$ on $2B$ for all
sufficiently large $k$. Caccioppoli's inequality therefore yields
\[
 \int_B|\nabla(w_k-u)|^2\,dm
 \lesssim r(B)^{-2}\int_{2B}|v_k-(v_k)_{2B}|^2\,dm\longrightarrow0.
\]
A finite covering of each compact subset of $\Omega$ by such balls
proves $\nabla w_k\to\nabla u$ strongly in $L^2_{\mathrm{loc}}(\Omega)$.


\section{Approximation for the rough Neumann problem}\label{sec.rough}

In this section we give the proof of Theorem \ref{thm.approxRoughNeumann}. 

\noindent\emph{Proof of Theorem \ref{thm.approxRoughNeumann}.}

Recall that for every $\xi\in\partial^*\Omega$,  the measure-theoretic outer normal $\nu(\xi)=(\nu_j(\xi))_{j=1}^{n+1}$ exists (see Section \ref{sec.geom}) and denote by $\sigma_*:=\mathcal H^n|_{\partial^* \Omega}=\sigma|_{\partial^* \Omega}$.

{\bf Part 1. Case $g$ is Lipschitz continuous.}

If $\Omega$ is a starlike Lipschitz domain, take $\Omega_1=\Omega$.
Otherwise, using Theorem \ref{thm.corona}, consider the countable sequence of starlike-Lipschitz subdomains $\{\Omega_i\}_i$, whose enumeration we fix. For $\sigma-$almost every $\xi\in\partial^*\Omega$, let $i_\xi$ be the index of the first
$\Omega_i$ in the enumeration such that $\xi\in\overline{\Omega_i}$, and define
\[
\hat\Gamma_\alpha(\xi):=\Gamma_\alpha^{\Omega_{i_\xi}}(\xi),\qquad\text{for }\xi\in\partial^*\Omega.
\]
It is easy to see that $\hat\Gamma_\alpha(\xi)\subseteq\Gamma_\alpha(\xi)$. On the set of $\sigma$-measure $0$ where $i_\xi$ is not well defined, simply set $\hat\Gamma_\alpha=\Gamma_\alpha^{\Omega}$.

\emph{Proof of \ref{item.c1}.} Suppose that $g\in\Lip(\partial^*\Omega)$ and let $\Phi=\Phi_g$ be the Varopoulos  extension of $g \, \chi_{\partial^* \Omega}$ onto $\Omega$ as in Lemma \ref{lm.varopoulos}. Now, for each $j\in\bb N$ large enough (namely, such that $2^{-j}<\diam(\Omega)/2$), let $\eta_j$ be a cut-off function satisfying that $\eta_j\equiv1$ on $B(\partial\Omega,2^{-(j+1)})$, $\supp\eta_j\subset B(\partial\Omega,2^{-j})$, $|\nabla\eta_j|\lesssim2^{j}$, and $0\leq\eta_j\leq1$. We set 
\begin{equation}\label{eq.fj0NeuRough}
\tilde H^{(k)}_j:=\rho^{(k)}_j\Phi\eta_j,\qquad\textrm{and}\qquad \tilde H^{(l)}_j:=\rho^{(l)}_j\Phi\eta_j
\end{equation}
where $\rho_j^{(k)}$ and $\rho_j^{(l)}$ are given by
\begin{align}\label{eq.rhoNeuRough} 
     \rho^{(k)}_j(x):=-\int_{\partial^*\Omega}\Big(\int_{\hat\Gamma_\alpha(\xi)}\eta_j\,dm\Big)^{-1}\chi_{\hat\Gamma_\alpha(\xi)}(x)\nu_l(\xi)d\sigma(\xi),&\qquad\textrm{and}
\\ \nonumber\rho^{(l)}_j(x):=-\int_{\partial^*\Omega}\Big(\int_{\hat\Gamma_\alpha(\xi)}\eta_j\,dm\Big)^{-1}\chi_{\hat\Gamma_\alpha(\xi)}(x)\nu_k(\xi)d\sigma(\xi).& 
\end{align}
In a manner analogous to  \eqref{eq:claimCBoundedByNNeu}, we claim that 
\begin{align}
\Vert\n C_\infty(\tilde H^{(k)}_j)\Vert_{L^p(\sigma)} + \Vert\n C_\infty(\tilde H^{(l)}_j)\Vert_{L^p(\sigma)}\lesssim\Vert \m M(\m N(\Phi))\Vert_{L^p(\sigma)}\lesssim\Vert g\Vert_{L^p(\sigma_*)}.\label{eq:claimCBoundedByNRoughNeu}
\end{align} 
The second estimate in this chain of estimates follows immediately by the construction of the Varopoulos extension. We postpone the proof of the first estimate in \eqref{eq:claimCBoundedByNRoughNeu} to the end of Part 1.

Given the claim (\ref{eq:claimCBoundedByNRoughNeu}), it follows that $\tilde H^{(k)}_j, \tilde H^{(l)}_j\in{\bf C}_{\infty,p}$, and by H\"older's inequality, in fact, we have that $\tilde H^{(k)}_j,\tilde H^{(l)}_j\in{\bf C}_{q,p}$ for all $q\in(1,\infty]$. Now, for each $q\in(1,\infty)$, by the density of $\Lip_c(\Omega)$ in ${\bf C}_{q,p}$ (see Lemma~\ref{lem:patch-density}), we may find compactly supported Lipschitz functions $H_j^{(k)},H_j^{(l)}$ such that 
\begin{equation}\label{eq.fjNeuRough}
\Vert\n C_q(H_j^{(k)}-\tilde H_j^{(k)})\Vert_{L^p(\sigma)}+\Vert\n C_q(H_j^{(l)}-\tilde H_j^{(l)})\Vert_{L^p(\sigma)}\leq2^{-j}\Vert g\Vert_{L^p(\sigma_*)}.
\end{equation}

This finishes the proof of \ref{item.c1} when $g$ is Lipschitz continuous.

\emph{Proof of \ref{item.c2}.} We continue to assume that  $g\in\Lip(\partial\Omega)$. Since $\Lip_c(\Omega)\subset L^2(\Omega)$, and $\partial_k H_j^{(k)}$ and $\partial_l H_j^{(l)}$ have zero mean-value, Lax-Milgram yields the existence of solutions $w_j\in \widehat{W}^{1,2}(\Omega)$ to \eqref{eq.poissonNeumann0Rough} with right hand side $\partial_kH^{(k)}_j-\partial_lH_j^{(l)}$.

Let $p'$ be the H\"older conjugate of $p$. By Proposition~\ref{prop.duality}, we know that $({\bf C}_{2,p'})^*={\bf N}_{2,p}$. Hence, we want to show that
\begin{equation}\label{eq.convNeu}
\int_\Omega (w_j-u)h\,dm\lra0,\quad\text{as }j\ra\infty, \text{ for each }h\in{\bf C}_{2,p'}.
\end{equation} 

Since both $w_j$ and $u$ have mean value zero on $\Omega$,  without loss of generality, we can restrict to $h\in{\bf C}_{2,p'}$ with mean value zero.

Assume first that $h\in\Lip_c(\Omega)$ with mean value zero, and let $V\in \widehat{W}^{1,2}(\Omega)$ be the weak solution to the Poisson problem
\begin{equation}\label{eq.poissonregularity}
	\left\{\begin{aligned}-\dv A^T\nabla V&=h,\qquad&\text{in }&\Omega,\\ \partial_{\nu_{A^T}}V&=0,\qquad&\text{on }&\partial\Omega.\end{aligned}\right. 
\end{equation}
It follows that
\begin{equation}\label{eq.conv2NeuRough}
	\int_\Omega (w_j-u)h\,dm=\int_\Omega A^T\nabla V\nabla(w_j-u)\,dm,
\end{equation}
where we used that $\partial_{\nu_{A^T}}V\equiv0$. On the other hand, note that
\begin{multline}\label{eq.conv3NeuRough}
	\int_\Omega A^T\nabla V\nabla(w_j-u)\,dm=\int_\Omega A^T\nabla V\nabla w_j\,dm-\int_\Omega A^T\nabla V\nabla u\,dm\\ =\int_\Omega A\nabla w_j\nabla V\,dm-\int_\Omega A\nabla u\nabla V\,dm=\int_\Omega H_j^{(l)} \partial_l V-H_j^{(k)} \partial_k V\,dm-(\partial_{\nu_A}u,V),
\end{multline}
since $Lu=0$, $\partial_{\nu_A}w_j\equiv 0$, $-\dv A\nabla w_j=  \partial_k H^{(k)}_j-\partial_l H_j^{(l)}$ in $\Omega$ and $V\in \widehat{W}^{1,2}(\Omega)$. Next, using \eqref{eq.fj0NeuRough}, and \eqref{eq.fjNeuRough}, we have
\begin{equation}\label{eq.conv31NeuRough}
\int_\Omega H_j^{(l)} \partial_l V-H_j^{(k)} \partial_k V\,dm=\int_\Omega (H_j^{(l)}-\tilde{H}_j^{(l)}) \partial_l V-(H_j^{(k)}-\tilde{H}_j^{(k)}) \partial_k V\,dm + \int_\Omega \tilde{H}_j^{(l)} \partial_l V-\tilde{H}_j^{(k)} \partial_k V\,dm.
\end{equation}

Putting (\ref{eq.conv2NeuRough}) together with  (\ref{eq.conv3NeuRough}), (\ref{eq.conv31NeuRough}), we get that
\begin{equation}\label{eq.conv4NeuRough}
	\int_\Omega (w_j-u)h\,dm=\int_\Omega (H_j^{(l)}-\tilde{H}_j^{(l)}) \partial_l V-(H_j^{(k)}-\tilde{H}_j^{(k)}) \partial_k V\,dm + \int_\Omega \tilde{H}_j^{(l)} \partial_l V-\tilde{H}_j^{(k)} \partial_k V\,dm~-~(\partial_{\nu_{A}}u,V).
\end{equation}
By  Carleson's Theorem (see Proposition \ref{prop.duality}), we obtain for the first term
\begin{align}
    &\Big|\int_\Omega (H_j^{(l)}-\tilde{H}_j^{(l)}) \partial_l V-(H_j^{(k)}-\tilde{H}_j^{(k)}) \partial_k V\,dm\Big|\nonumber
    \\
    &\lesssim\big(\Vert\n C_2(H_j^{(k)}-\tilde H_j^{(k)})\Vert_{L^p(\sigma)} + \Vert\n C_2(H_j^{(k)}-\tilde H_j^{(k)})\Vert_{L^p(\sigma)}\big)\Vert\widetilde{\m N}_2(\nabla V)\Vert_{L^{p'}(\sigma)} \nonumber\\
    &    \lesssim 2^{-j}\Vert g\Vert_{L^p(\sigma_*)}\Vert\n C_2(h)\Vert_{L^{p'}(\sigma)} \label{eq.fjapproxNeuRough},
\end{align}
where we used that the Poisson-Neumann-Regularity problem $(\operatorname{PNR}_{p'}^{L^*})$ is solvable in $\Omega$, since $(\operatorname{N}_p^L)$ and $(\operatorname{D}_{p'}^{L^*})$ are solvable in $\Omega$ (see \cite[Definition 1.13, and Lemma 7.1]{fl24}).

The second and third term of \eqref{eq.conv4NeuRough} can be combined. First, we use the definition of $\partial_{\nu_A}u\equiv\partial_{t,l,k}g$ to split  
\[-(\partial_{\nu_A}u,V )=\int_{\partial^*\Omega}g\partial_{t,l,k}V\, d\sigma=\int_{\partial^*\Omega}g\nu_l\partial_kV\,d\sigma-\int_{\partial^*\Omega}g\nu_k\partial_lV\, d\sigma=:I_k+I_l.\]
Next, we split
\[\int_\Omega \tilde{H}_j^{(l)} \partial_l V-\tilde{H}_j^{(k)} \partial_k V\,dm=\int_\Omega \tilde{H}_j^{(l)} \partial_l V\,dm-\int_\Omega\tilde{H}_j^{(k)} \partial_k V\,dm=:II_l+II_k.\]

Now, we combine $I_k+II_k$ and $I_l+II_l$, but since the estimates work completely analogously just with a changed sign, we only demonstrate how to bound $I_k+II_k$. To that extent, we use the definition of $\rho_j^{(k)}$ (see \eqref{eq.rhoNeuRough}) and Fubini-Tonelli to rewrite $II_k$ as
\begin{align*}
    II_k &=
    -\int_\Omega \rho_j^{(k)}\Phi\eta_j \partial_kv\,dm
    \\
    &=-\int_{\Omega} \left[\int_{\partial^*\Omega}\Big(\int_{\hat\Gamma_\alpha(\xi)}\eta_j\,dm\Big)^{-1}\chi_{\hat\Gamma_\alpha(\xi)}(x)\nu_l(\xi)\,d\sigma(\xi)\right]\Phi(x)\eta_j(x) \partial_k V(x)\,dm(x)
    \\
    &=-\int_{\partial^*\Omega}\left[\int_\Omega\Big(\int_{\hat\Gamma_\alpha(\xi)}\eta_j\,dm\Big)^{-1}\chi_{\hat\Gamma_\alpha(\xi)}(x)\Phi(x)\eta_j(x) \partial_k V(x)\,dm(x)\right]\nu_l(\xi)\,d\sigma(\xi)
    \\
    &=-\int_{\partial^*\Omega}\left[\int_{\hat\Gamma_\alpha(\xi)}\Big(\int_{\hat\Gamma_\alpha(\xi)}\eta_j\,dm\Big)^{-1}\Phi(x)\eta_j(x) \partial_k V(x)\,dm(x)\right]\nu_l(\xi)\,d\sigma(\xi).
\end{align*}
As in the proof of Theorem~\ref{thm.approxNeumann}, for $x \in \Omega$ and $\xi \in \partial^*\Omega$, we set \
$$
m_j(x,\xi):=\Big(\int_{\hat\Gamma_\alpha(\xi)}\eta_j\,dm\Big)^{-1}\eta_j(x),
$$
which satisfies \(\supp(m_j(\cdot,\xi)\subset B(\partial\Omega,2^{-j})\) and \(\int_{\hat\Gamma_\alpha(\xi)}m_j(x,\xi)dx=1\). Combining $II_k$ with $I_k$ gives
\begin{align*}
    I_k &+II_k 
    =\int_{\partial^*\Omega} g(\xi)\partial_kV(\xi)\nu_l(\xi)\,d\sigma(\xi)-\int_{\partial^*\Omega}\left[\int_{\hat\Gamma_\alpha(\xi)}m_j(x,\xi)\Phi(x)\partial_kV(x)\,dm(x)\right]\nu_l(\xi)\,d\sigma(\xi)
    \\
    &=\int_{\partial^*\Omega}\left[\int_{\hat\Gamma_\alpha(\xi)}m_j(x,\xi)\big(g(\xi)\partial_kV(\xi)-\Phi(x) \partial_kV(x)\big)\,dm(x)\right]\nu_l(\xi)\,d\sigma(\xi)
    \\
    &=\int_{\partial^*\Omega}\left[\int_{\hat\Gamma_\alpha(\xi)}m_j(x,\xi)\big[g(\xi)(\partial_kV(\xi)-\partial_kV(x))+(g(\xi)-\Phi(x)) \partial_kV(x)\big]\,dm(x)\right]\nu_l(\xi)\,d\sigma(\xi).\label{eq.conv5NeuRough}
\end{align*}
For $I_l+II_l$ we can proceed analogously. This gives us for the sum of the second and third term of \eqref{eq.conv4NeuRough}
\begin{align*}
   I_k &+II_k+I_l+II_l
    \\
    &= \int_{\partial^*\Omega}\left[\int_{\hat\Gamma_\alpha(\xi)}m_j(x,\xi)\big[g(\xi)(\partial_kV(\xi)-\partial_kV(x))+(g(\xi)-\Phi(x)) \partial_k V(x)\big]\,dm(x)\right]\nu_l(\xi)\,d\sigma(\xi)
    \\
    &-\int_{\partial^*\Omega}\left[\int_{\hat\Gamma_\alpha(\xi)}m_j(x,\xi)\big[g(\xi)(\partial_lV(\xi)-\partial_lV(x))+(g(\xi)-\Phi(x)) \partial_l V(x)\big]\,dm(x)\right]\nu_k(\xi)\,d\sigma(\xi).
\end{align*}
Reorganizing the terms, gives
\begin{align*}
     &I_k+II_k+I_l+II_l
    \\
    &= \int_{\partial^*\Omega}\left[\int_{\hat\Gamma_\alpha(\xi)}m_j(x,\xi)\big[g(\xi)(\partial_kV(\xi)-\partial_kV(x))\nu_l(\xi)-g(\xi)(\partial_lV(\xi)-\partial_lV(x))\nu_k(\xi)\big]\,dm(x)\right]\,d\sigma(\xi)
    \\
    &\qquad+\int_{\partial^*\Omega}\left[\int_{\hat\Gamma_\alpha(\xi)}m_j(x,\xi)\big[(g(\xi)-\Phi(x)) \partial_k V(x)\nu_l(\xi)-(g(\xi)-\Phi(x)) \partial_l V(x)\nu_k(\xi)\big]\,dm(x)\right]\,d\sigma(\xi)
    \\
    &=:J_1+J_2.
\end{align*}
For $J_2$, we use that \(\Phi\) is Lipschitz as the Varopolous extension of $g$ (see Lemma \ref{lm.varopoulos}). Hence, we have $|g(\xi)-\Phi(x)|\leq L2^{-j}$, where $L$ is the Lipschitz constant of $\Phi$, and 
\begin{equation}
    |J_2|\leq L2^{-j}\Vert \widehat{\mathcal{N}}^{2^{-j+2}}(\nabla V)\Vert_{L^1(\sigma)},
\end{equation}
where $\widehat{\mathcal{N}}^{2^{-j+2}}(\nabla V)(\xi):=\sup_{y\in\Gamma_\alpha(\xi)\cap B(\xi, 2^{-j+2})}\fint_{B(y,\delta(y)/2)}|\nabla V(y)|$ is the truncated $L^1-$mean-valued nontangential maximal function. Since $\widetilde{\mathcal{N}}^{2^{-j+2}}(\nabla V)(\xi)\leq \widetilde{\mathcal{N}}_{2}(\nabla V)(\xi)$, we can use the solvability of the Poisson-Neumann-regularity problem $(\operatorname{PNR})_{p'}^{L^*}$, and obtain
\begin{align}
    |J_2|\leq L2^{-j}\Vert \widehat{\mathcal{N}}^{2^{-j+2}}(\nabla V)\Vert_{L^1(\sigma)}\lesssim 2^{-j}\Vert \widehat{\mathcal{N}}_{2}(\nabla V)\Vert_{L^{p'}(\sigma)}\lesssim 2^{-j}\Vert \n C_2(h)\Vert_{L^{p'}(\sigma)}.
\end{align}

For $J_1$, we note that for every $\xi \in \partial^* \Omega$ and every $x\in \mathrm{supp}(\eta_j)\cap\Gamma_\alpha(\xi)$, we have that
\begin{multline*}
    g(\xi)(\partial_kV(\xi)-\partial_kV(x))\nu_l(\xi)-g(\xi)(\partial_lV(\xi)-\partial_lV(x))\nu_k(\xi)
    \\
    = g(\xi)\Big[\big(\partial_kV(\xi)\nu_l(\xi)-\partial_lV(\xi)\nu_k(\xi)\big) -\big(\partial_kV(x)\nu_l(\xi)-\partial_lV(x)\nu_k(\xi)\big)\Big]\\=g(\xi)\big[\partial_{t,l,k}V(\xi)-\nabla V(x)\cdot \vec{T}(\xi) \big],
\end{multline*}
where 
\[
\vec{T}_{l,k}(\xi):=\nu_l(\xi)\vec{e}_k-\nu_{k}(\xi)\vec{e}_l
\]
is a tangent vector to $\partial\Omega$ at $\xi\in \partial^*\Omega$.


Now, suppose momentarily that $\Omega$ is a starlike Lipschitz domain. In that case, $\pom = \partial^*\Omega$ and  \cite[Theorem 3.1]{kp}   yields that $\fint_{B(y_\xi,\delta(y_\xi))}\nabla V(x)\cdot \vec{T}_{l,k}(\xi)\,dm(x)$ with $y_\xi\in \Gamma_\alpha(\xi)$  converges weakly in $L^{p'}(\partial\Omega)$ to $\nabla V(\xi)\cdot \vec{T}_{l,k}(\xi)=\partial_{t,l,k}V$ on $\partial\Omega$ as $y_\xi\ra\xi$.  Since $g\in L^p$, using the definition of $\eta_j$, it easily follows that  
\[J_1=-\int_{\partial\Omega}\left[\int_{\Gamma_\alpha(\xi)}\Big(\int_{\Gamma_\alpha(\xi)}\eta_j\,dm\Big)^{-1}\eta_j(x)(\partial_{t,l,k}V(\xi)-\nabla V(x)\cdot \vec{T}_{l,k}(\xi))\,dm(x)\right]g(\xi)\,d\sigma(\xi)\to 0\]
as $j\to \infty$.   Indeed, fix $\varepsilon>0$ small enough to be chosen and split 
$$\Gamma_\alpha(\xi)=\big[\Gamma_\alpha(\xi) \cap B(\xi,\varepsilon 2^{-j})\big] \bigcup \big[\Gamma_\alpha(\xi) \setminus B(\xi,\varepsilon 2^{-j})\big]=: \gamma^\ve(\xi) \cup w_\ve(\xi),$$
which introduces the corresponding splitting $J_1= J_{11}(\ve)+J_{12}(\ve).$  Using the fact that $m(\gamma^\ve(\xi))\lesssim (\ve 2^{-j})^{n+1}$, we infer that
$$
|J_{11}(\ve)| \lesssim \ve^{n+1} \left[ \|\partial_{t,l,k}V\|_{L^{p'}(\partial \Omega)} + \|\widetilde{\mathcal{N}}_2(\nabla V)\|_{L^{p'}(\partial \Omega)} \right] \|g\|_{L^p(\partial \Omega)}.
$$
Also, noting that $w_\ve(\xi)$ is  a Whitney region, by the weak convergence in $L^{p'}(\partial \Omega)$, we have that $J_{12}(\ve) \to 0$ as $j \to \infty$. Therefore, 
$$
|J_1| \lesssim \ve^{n+1} \left[ \|\partial_{t,l,k}V\|_{L^{p'}(\partial \Omega)} + \|\widetilde{\mathcal{N}}_2(\nabla V)\|_{L^{p'}(\partial \Omega)} \right] \|g\|_{L^p(\partial \Omega)} + |J_{12}(\ve)|$$
and by taking limits first as $j \to \infty$ and then as $\ve \to 0$, we conclude the desired convergence.

In total, combining this with \eqref{eq.conv4NeuRough} and \eqref{eq.fjapproxNeuRough} yields
\[\int_\Omega (w_j-u)h\,dm\lesssim 2^{-j}\Vert g\Vert_{L^p(\sigma)} + |J_1| + |J_2|\lesssim 2^{-j}(1+\Vert g\Vert_{L^p(\sigma)}) + |J_1|\to 0\]
as $j \to \infty$. This concludes the proof of Part 1 in the case that $\Omega$ is a starlike Lipschitz domain.
 
If $\Omega$ is not a starlike Lipschitz domain, then as mentioned at the beginning of the proof, we may apply Theorem \ref{thm.corona} and consider the sequence of starlike-Lipschitz subdomains $\{\Omega_i\}_i$, whose enumeration has been fixed, with the property that $\partial^*\Omega$ can be exhausted by the union of the $\partial\Omega_i$'s up to a set of $\sigma$ measure $0$. Recall that $i_\xi$ is defined as the   index of the first $\Omega_i$ for which $\xi\in\overline{\Omega_i}$. Let $E_i\subset\partial^*\Omega$ consist  of the points $\xi\in\partial^*\Omega$ such that $i_\xi=i$. Thus we may write $E_i$ as
\[
E_i=(\partial^*\Omega\cap\partial\Omega_i)\setminus\cup_{m<i}\partial\Omega_m,
\]
Accordingly, given $I\in\bb N$ large, we may decompose $J_1$ as follows:


\begin{equation}\label{eq.decomp1}
\begin{split}
|J_1|&\leq\sum_{i=1}^I\Big|\int_{E_i}\Big(\int_{\hat\Gamma_\alpha(\xi)}m_j(x)\big[\partial_{t,l,k}V(\xi)-\nabla V(x)\cdot \vec{T}_{l,k}(\xi) \big]dx\Big)g(\xi)\,d\sigma(\xi)\Big|
\\
~&\qquad+~\Big|\int_{\partial^*\Omega\setminus(\cup_{i=1}^I\partial\Omega_i)}\Big(\int_{\hat\Gamma_\alpha(\xi)}m_j(x)\big[\partial_{t,l,k}V(\xi)-\nabla V(x)\cdot \vec{T}_{l,k}(\xi) \big]dx\Big)g(\xi)\,d\sigma(\xi)\Big|=:J_{11}+J_{12}.
\end{split}
\end{equation}

Note that the integrand in $J_{12}$ is pointwise controlled from above by an $L^1(\partial^*\Omega)$ function independent of $i$ and $j$. Thus, by the Lebesgue Dominated Convergence Theorem, we may choose $I$ depending on $j$ so large that $\sigma\big(\partial^*\Omega\setminus(\cup_{i=I}^\infty\partial\Omega_i)\big)$ is sufficiently small. Having fixed $I$ so that $J_{12}$ is as small as desired, we turn to $J_{11}$, which is a finite sum of integrals.

Now, let $j$ be large enough so that $2^{2-j} < \min_{i=1, \dots, I} \diam\Omega_i$. By the definition of $E_i$ and $\hat\Gamma_\alpha(\xi)$, we have that
\begin{equation}\label{eq.decomp2}\nonumber
\begin{split}
&\Big|\int_{E_i}\Big(\int_{\hat\Gamma_\alpha(\xi)}m_j(x)\big[\partial_{t,l,k}V(\xi)-\nabla V(x)\cdot \vec{T}_{l,k}(\xi) \big]dx\Big)g(\xi)\,d\sigma(\xi)\Big|
\\
&=\Big|\int_{(\partial^*\Omega\cap\partial\Omega_i)\setminus\cup_{m<i}\partial\Omega_m}\Big(\int_{\Gamma^{\Omega_i}_\alpha(\xi)}m_j(x)\big[\partial_{t,l,k}V(\xi)-\nabla V(x)\cdot \vec{T}_{l,k}(\xi) \big]dx\Big)g(\xi)\,d\sigma(\xi)\Big|.
\end{split}
\end{equation}
We now use that $\Omega_i$ is a starlike Lipschitz subdomain, and proceed as before to see that $\lim_{j\ra\infty}J_{11}=0$. Here the only change needed is the observation that \cite[Theorem 3.1]{kp} can also be applied on only the part \(E_i\) of the boundary for which we have \(L^p\) control of \(\tilde{N}(\nabla u)\). Since there are finitely many $i$, the desired result is attained, modulo the proof of Claim \eqref{eq:claimCBoundedByNRoughNeu}.

It remains to prove the claim \eqref{eq:claimCBoundedByNRoughNeu}. We only provide the proof for \(\tilde H_j^{(k)}\) since the proof for \(\tilde H_j^{(l)}\) is analogous. For $\sigma$-almost every $\xi\in\partial^*\Omega$, set	$\tau_{j,\xi}=\min\{2^{-j}, \diam\Omega_{i_\xi}\}$. It is easy to see that
	\[
	\int_{\hat\Gamma_\alpha(\xi)}\eta_j\,dm\approx\tau_{j,\xi}^{n+1}.
	\]
Hence, by \eqref{eq.rhoNeuRough}, we have
\begin{align*}
\sup_{y\in B(x,\hat c \delta(x))}|\tilde{H}_j^{(k)}(y)|&=\sup_{y\in B(x,\hat c \delta(x))}
\left(
\int_{\partial^*\Omega}\Big(\int_{\hat\Gamma_\alpha(\xi)}\eta_j\,dm\Big)^{-1}
\chi_{\hat\Gamma_\alpha(\xi)}(y)\eta_j(y)|\Phi(y)|  
\,d\sigma(\xi)
\right)
\\
&\leq 
\int_{\partial^*\Omega}
\tau_{j,\xi}^{-(n+1)}\big(\sup_{y\in B(x,\hat c \delta(x))}\chi_{\hat\Gamma_\alpha(\xi)}(y)\eta_j(y)\big)\sup_{y\in \hat\Gamma_\alpha(\xi)}|\Phi(y)| 
\,d\sigma(\xi).
\\
&\leq 
\int_{\partial^*\Omega}
\tau_{j,\xi}^{-(n+1)}\big(\sup_{y\in B(x,\hat c \delta(x))}\chi_{\hat\Gamma_\alpha(\xi)}(y)\eta_j(y)\big)\mathcal{N}(\Phi)(\xi) 
\,d\sigma(\xi).
\end{align*}
Here we used \(\hat\Gamma_\alpha(\xi)\subset\Gamma_\alpha(\xi)\) to get the nontangential maximal function.
Hence, for a fix $\zeta\in\partial\Omega$, we have
\begin{align}
\n C_\infty(\tilde H_j^{(a)})(\zeta)
&\leq \sup_{r>0}\frac{1}{r^n}
\int_{B(\zeta,r)\cap\Omega}
\sup_{y\in B(x,\hat c \delta(x))}|\rho_j(y)|\,dm(x)\nonumber
\\
&\leq \sup_{r>0}\frac{1}{r^n}
\int_{B(\zeta,r)\cap\Omega}
\Big[\int_{\partial^*\Omega}
\tau_{j,\xi}^{-(n+1)}\Big(\sup_{y\in B(x,\hat c \delta(x))}\chi_{\hat\Gamma_\alpha(\xi)}(y)\eta_j(y)\Big)\mathcal{N}(\Phi)(\xi) 
\,d\sigma(\xi)\Big]\,dm(x).\label{eq.HjRoughNeuB}
\end{align}
By Fubini-Tonelli, we continue \eqref{eq.HjRoughNeuB} with
\begin{align*}
\n C_\infty(\tilde H_j^{(a)})(\zeta)
&\leq \sup_{r>0}\frac{1}{r^n}
\Big[\int_{\partial^*\Omega}
\tau_{j,\xi}^{-(n+1)}\Big(\int_{B(\zeta,r)\cap\Omega}\sup_{y\in B(x,\hat c \delta(x))}\chi_{\hat\Gamma_\alpha(\xi)}(y)\eta_j(y)\,dm(x)\Big)\mathcal{N}(\Phi)(\xi) 
\,d\sigma(\xi)\Big].
\end{align*}
We can note that there exists a constant \(C_\alpha>0\) such that if \(\xi\notin B(\zeta,C_\alpha r)\cap\partial\Omega\) and \(y\in B(\zeta,2r)\cap\Omega\), then \(\chi_{\hat\Gamma_\alpha(\xi)}(y)=0\). Since also 
\[\int_{B(\zeta,r)\cap\Omega}\sup_{y\in B(x,\hat c \delta(x))}\chi_{\hat\Gamma_\alpha(\xi)}(y)\eta_j(y)\,dm(x)\lesssim \int_{\hat\Gamma_\alpha(\xi)}\eta_j(y)\,dm(x)\lesssim \tau_{j,\xi}^{n+1},\]
we obtain
\begin{align*}
\n C_\infty(\tilde H_j^{(a)})(\zeta)
&\leq \sup_{r>0}\frac{1}{r^n}
\Big[\int_{B(\zeta,C_\alpha r)\cap\partial\Omega}
\tau_{j,\xi}^{-(n+1)} \tau_{j,\xi}^{n+1}\mathcal{N}(\Phi)(\xi) 
\,d\sigma(\xi)\Big]
\lesssim \mathcal{M}(\mathcal{N}(\Phi))(\zeta).
\end{align*}
This completes the proof of claim \eqref{eq:claimCBoundedByNRoughNeu}.

{\bf Part 2. Case $g\in L^p(\partial^*\Omega,\sigma)$.}


Now, we fix $g\in L^p(\partial^*\Omega,\sigma)$, and let $j\in\mathbb{N}$. Since $\Lip(\partial^*\Omega,\sigma)$ is dense in $L^p(\partial^*\Omega,\sigma)$, there exists $g_j\in\Lip(\partial^*\Omega)$ with $\Vert g_j-g\Vert_{L^p}<\frac{1}{j}$. Let us denote by $u_j$ the solution to \eqref{eq.neumann0Rough} with boundary data $g_j$. Let \(\{h_m\}_{m\in\mathbb{N}}\subset\Lip_c(\Omega)\) be a countable subset that is dense in \({\bf C}_{2,p'}(\Omega)\) (see Lemma~\ref{lem:patch-separability}). Furthermore by Part 1, there exists $\tilde{H}_j^{(k), i},\tilde{H}_j^{(l), i}\in\Lip_c(\Omega)$ and solutions $w_j^i\in \widehat{W}^{1,2}(\Omega)$ to \eqref{eq.poissonNeumann0Rough} such that $\Vert \tilde{H}_j^{(k),i}\Vert_{{\bf C}_{2,p'}}+\Vert \tilde{H}_j^{(l),i}\Vert_{{\bf C}_{2,p'}}\leq C\Vert g_j\Vert_{L^p(\sigma_*)}$ and $\tilde{w}_j^i\ra  u_j$ as $i\ra\infty$ in the weak-$*$ topology in ${\bf N}_{2,p}$. For every \(j\), we can now choose \(i_j\) such that 
\begin{equation}\label{eq.flargRough}
\int_\Omega (\nabla \tilde{w}_j^{i_j}-\nabla u_j)F_m \,dm\leq \frac{1}{j} \qquad \textrm{  for all }1\leq m\leq j.
\end{equation}
We set $H_j^{(k)}:=\tilde{H}_j^{(k), i_j}$, $H_j^{(l)}:=\tilde{H}_j^{(l), i_j}$ and \(w_j:=\tilde{w}_j^{i_j}\), and we claim that the sequences $H_j^{(k)}, H_j^{(l)}$ and $w_j$ satisfy the claims of the theorem. Indeed,
\[\Vert H_j^{(k)}\Vert_{{\bf C}_{2,p'}}+\Vert H_j^{(l)}\Vert_{{\bf C}_{2,p'}}=\Vert \tilde{H}_j^{(k), i_j}\Vert_{{\bf C}_{2,p'}}+\Vert \tilde{H}_j^{(l),i_j}\Vert_{{\bf C}_{2,p'}}\leq C\Vert g_j\Vert_{L^p}\leq C(\Vert g-g_j\Vert_{L^p}+\Vert g\Vert_{L^p})\leq 2C\Vert g\Vert_{L^p},\]
which shows \eqref{eq.uniformc}. Moreover, for every \(\varepsilon>0\) and for every \(h\in {\bf C}_{2,p'}\) there exists \(j_0>\frac{1}{\varepsilon}\) such that \(\Vert h_{j_0}-h\Vert_{{\bf C}_{2,p'}}<\varepsilon\). Then,  for all \(j\geq j_0\), we have
\begin{align*}
    \int_\Omega (w_j-u)h\, dm&=\int_\Omega (w_j-u)(h-h_{j_0})\, dm + \int_\Omega (w_j-u)h_{j_0}\, dm=I_1+I_2.
\end{align*}
By \eqref{eq.flargRough}, Proposition~\ref{prop.duality}, and the solvability  of the rough Neumann problem, we have
\[|I_2|\leq \Big|\int_\Omega (w_j-u_j)h_{j_0}\, dm\Big| + \Big|\int_\Omega (u_j-u)h_{j_0}\, dm\Big|\leq \frac{1}{j}+\Vert u_j-u\Vert_{{\bf N}_{2,p}}\leq \frac{1}{j}+\Vert g_j-g\Vert_{L^p}\leq \frac{2}{j_0}<2\varepsilon,\]
and, by Proposition~\ref{prop.duality} and the solvability of \((\N_p^L)\) and of the Poisson-Neumann problem $(\operatorname{PN}_p^L)$, which follows from the solvability of $(\N_{p'}^{L^*})$ and $(\Di_{p}^{L})$ (see \cite[Lemma 7.1]{fl24}), we get
\begin{align*}
    |I_1|&\leq \Big|\int_\Omega (w_j-u)(h-h_{j_0})\, dm\Big|
    \\
    &\leq(\Vert w_j\Vert_{{\bf N}_{2,p}}+\Vert u\Vert_{{\bf N}_{2,p}})\Vert h_{j_0}-h\Vert_{{\bf C}_{2,p'}}\lesssim \varepsilon(\Vert H_j\Vert_{{\bf C}_{2,p}}+\Vert g\Vert_{L^p})\lesssim \varepsilon\Vert g\Vert_{L^p}.
\end{align*}
This completes the proof of \(w_j\to u\) in \({\bf N}_{2,p}\).

Fix a ball $B$ with $\overline{4B}\subset\Omega$. By
construction of \(H_j^{(k)}\) and \(H_j^{(l)}\), $Lw_j=0$ on $4B$ for all
sufficiently large $j$, and $\{w_j\}_j$ is uniformly bounded in
$L^2(4B)$. Caccioppoli's inequality gives a uniform
$W^{1,2}(2B)$ bound on \(\{w_j\}_j\). The Rellich–Kondrachov theorem then gives relative compactness in \(L^2(2B)\). Every subsequential limit equals \(u\), since the weak-* convergence in \({\bf N}_{2,p}\) implies
\[\int_{2B}w_j \phi \,dm\to \int_{2B}u \phi \,dm\qquad \textrm{ for all }\phi\in C^\infty_c(2B).\] Consequently,
\[
 w_j\longrightarrow u\quad\text{strongly in }L^2(2B).
\]
Since \(j\) was chosen sufficiently large such that \(L(w_j-u)=0\) on \(4B\), we get by Caccioppoli
\[
 \int_B|\nabla(w_j-u)|^2\,dm
       \lesssim r(B)^{-2}\int_{2B}|w_j-u|^2\,dm\longrightarrow0,\qquad \textrm{as j}\to \infty,
\]
and by De Giorgi--Nash--Moser theory for some $\alpha>0$
\[
 \|w_j-u\|_{C^\alpha(\overline B)}
       \le C_B\|w_j-u\|_{L^2(2B)}\longrightarrow0.
\]
Covering compact subsets of $\Omega$ by finitely many such balls
proves the strong convergence in $W^{1,2}_{\mathrm{loc}}(\Omega)$
and in $C^\alpha_{\mathrm{loc}}(\Omega)$.

\appendix

\section{Auxiliary results}

\begin{lemma}\label{lem:patch-density}
Let $q,p\in(1,\infty)$, $H\in\mathbf C_{q,p}$, and $\varepsilon>0$.
There exists $H_\varepsilon\in\operatorname{Lip}_c(\Omega)$ such that
\[
 \|H-H_\varepsilon\|_{\mathbf C_{q,p}}<\varepsilon.
\]
If $\int_\Omega H\,dm=0$, then $H_\varepsilon$ may also be chosen to
have integral zero. Moreover, if $H=0$ almost everywhere on
$\{x\in\Omega:\delta_\Omega(x)\ge r\}$, then $H_\varepsilon$ may be chosen so that
\[
 \operatorname{supp}H_\varepsilon
 \Subset\{x\in\Omega:\delta_\Omega(x)<2r\},
\]
with the integral-zero condition preserved when applicable.
\end{lemma}

\begin{proof}
The first assertion is \cite[Lemma 2.6]{mpt25}, the density of
$\operatorname{Lip}_c(\Omega)$ in $\mathbf C_{q,p}$.
To obtain the support assertion, choose a Lipschitz function
$\chi:\Omega\to[0,1]$ that equals one on $\{\delta_\Omega\le r\}$
and vanishes on $\{\delta_\Omega\ge 3r/2\}$.
If no support restriction is imposed, set $\chi\equiv1$.
Given $\tau>0$, choose $K\in\operatorname{Lip}_c(\Omega)$ with
$\|H-K\|_{\mathbf C_{q,p}}<\tau$, and set $\widetilde H=\chi K$.
Since $\chi H=H$ and multiplication by a function bounded by one
does not increase the Carleson norm,
\[
 \|H-\widetilde H\|_{\mathbf C_{q,p}}<\tau.
\]
This proves the support assertion when no integral condition is imposed.

Suppose now that $\int_\Omega H\,dm=0$. Fix
$\varphi\in\operatorname{Lip}_c(\Omega)$ with
$\int_\Omega\varphi\,dm=1$; when the support restriction is imposed,
choose $\operatorname{supp}\varphi\Subset\{\delta_\Omega<2r\}$.
Proposition~\ref{prop.duality} applied to $\mathbf 1_\Omega$,
gives
\[
 \left|\int_\Omega\widetilde H\,dm\right|
 =\left|\int_\Omega(\widetilde H-H)\,dm\right|
 \le C_0\|\widetilde H-H\|_{\mathbf C_{q,p}}<C_0\tau,
\]
where $C_0$ is finite because $\Omega$ is bounded. Set
\[
 H_\varepsilon
 =\widetilde H-
       \left(\int_\Omega\widetilde H\,dm\right)\varphi.
\]
Then $H_\varepsilon$ has integral zero, has the required support, and
\[
 \|H-H_\varepsilon\|_{\mathbf C_{q,p}}
 <\bigl(1+C_0\|\varphi\|_{\mathbf C_{q,p}}\bigr)\tau.
\]
Choose $\tau<\varepsilon/
(1+C_0\|\varphi\|_{\mathbf C_{q,p}})$.
\end{proof}

\begin{lemma}\label{lem:patch-separability}
For $q,p\in(1,\infty)$, the unit ball of $\mathbf C_{q,p}$ has a
countable dense subset contained in $\operatorname{Lip}_c(\Omega)$.
The same assertion holds for any finite Cartesian power of this space.
\end{lemma}

\begin{proof}
Since $\Omega$ is bounded,
\[
 \|H\|_{\mathbf C_{q,p}}
 \le C_\Omega\|H\|_{L^\infty(\Omega)}
 \qquad(H\in C_0(\Omega)).
\]
Indeed, $\mathscr C_q(H)\le C\operatorname{diam}(\Omega)
\|H\|_\infty$ pointwise on $\partial\Omega$.
The space $C_0(\Omega)$ is separable in the uniform norm, and
$\operatorname{Lip}_c(\Omega)$ is dense in it. Thus there is a
countable subset of $\operatorname{Lip}_c(\Omega)$ that is dense in
$\mathbf C_{q,p}$. Applying the continuous map
$H\mapsto H/\max\{1,\|H\|_{\mathbf C_{q,p}}\}$ gives a countable
dense subset of its unit ball. Finite products are handled in the same way.
\end{proof}

\bibliographystyle{alpha-sort-max} 
\bibliography{refs} 

@article {mz25,
    AUTHOR = {Mourgoglou, Mihalis and Zacharopoulos, Thanasis},
     TITLE = {Varopoulos extensions in domains with {A}hlfors-regular
              boundaries and applications to boundary value problems for
              elliptic systems with {$L^\infty$} coefficients},
   JOURNAL = {Adv. Math.},
  FJOURNAL = {Advances in Mathematics},
    VOLUME = {461},
      YEAR = {2025},
     PAGES = {Paper No. 110054, 85},
      ISSN = {0001-8708,1090-2082},
   MRCLASS = {42B35 (35J47)},
  MRNUMBER = {4833689},
MRREVIEWER = {Lubomira\ G.\ Softova},
       DOI = {10.1016/j.aim.2024.110054},
       URL = {https://doi.org/10.1016/j.aim.2024.110054},
}

@article {haj96,
    AUTHOR = {Haj\l asz, Piotr},
     TITLE = {Sobolev spaces on an arbitrary metric space},
   JOURNAL = {Potential Anal.},
  FJOURNAL = {Potential Analysis. An International Journal Devoted to the
              Interactions between Potential Theory, Probability Theory,
              Geometry and Functional Analysis},
    VOLUME = {5},
      YEAR = {1996},
    NUMBER = {4},
     PAGES = {403--415},
      ISSN = {0926-2601},
   MRCLASS = {46E35},
  MRNUMBER = {1401074},
MRREVIEWER = {J. Musielak},
       DOI = {10.1007/BF00275475},
       URL = {https://doi.org/10.1007/BF00275475},
}

@Article{ht21,
  author	= {Hofmann, Steve and Tapiola, Olli},
  title		= {Uniform rectifiability implies {V}aropoulos extensions},
  journal	= {Adv. Math.},
  fjournal	= {Advances in Mathematics},
  volume	= {390},
  year		= {2021},
  pages		= {Paper No. 107961, 53},
  issn		= {0001-8708},
  mrclass	= {42B37 (31B15 31E05)},
  mrnumber	= {4298596},
  mrreviewer	= {Tuomas P. Hyt\"{o}nen},
  doi		= {10.1016/j.aim.2021.107961},
  url		= {https://doi.org/10.1016/j.aim.2021.107961}
}

@article {hr13,
    AUTHOR = {Hyt\"{o}nen, Tuomas and Ros\'{e}n, Andreas},
     TITLE = {On the {C}arleson duality},
   JOURNAL = {Ark. Mat.},
  FJOURNAL = {Arkiv f\"{o}r Matematik},
    VOLUME = {51},
      YEAR = {2013},
    NUMBER = {2},
     PAGES = {293--313},
      ISSN = {0004-2080},
   MRCLASS = {42B37 (35J25)},
  MRNUMBER = {3090198},
MRREVIEWER = {Michael A. Perelmuter},
       DOI = {10.1007/s11512-012-0167-7},
       URL = {https://doi.org/10.1007/s11512-012-0167-7},
}

@article {hr18,
    AUTHOR = {Hyt\"onen, Tuomas and Ros\'en, Andreas},
     TITLE = {Bounded variation approximation of {$L_p$} dyadic martingales
              and solutions to elliptic equations},
   JOURNAL = {J. Eur. Math. Soc. (JEMS)},
  FJOURNAL = {Journal of the European Mathematical Society (JEMS)},
    VOLUME = {20},
      YEAR = {2018},
    NUMBER = {8},
     PAGES = {1819--1850},
      ISSN = {1435-9855,1435-9863},
   MRCLASS = {42B37 (35J15 42B25)},
  MRNUMBER = {3854892},
MRREVIEWER = {Giuseppe\ Di Fazio},
       DOI = {10.4171/JEMS/800},
       URL = {https://doi.org/10.4171/JEMS/800},
}

@unpublished {fl24,
 AUTHOR = {Feneuil, Joseph and Li, Linhan},
     TITLE = {The ${L}^p$ {P}oisson-{N}eumann problem and its relation to the {N}eumann problem},
      NOTE ={Preprint. June 2024. arXiv: 2406.16735},
     URL = {https://arxiv.org/abs/2406.16735},
}

@unpublished {fyy25,
 AUTHOR = {Fu, Xuelian and Yang, Dachun and Yang, Sibei},
     TITLE = {Equivalent {C}haracterizations and {T}heir {A}pplications of {S}olvability of ${L}^p$ {P}oisson--{R}obin(-{R}egularity) {P}roblems on Rough Domains},
      NOTE ={Preprint. July 2025. arXiv: 2507.11103},
     URL = {https://arxiv.org/abs/2507.11103},
}

@misc{bmp25,
      title={The {P}oisson-{D}irichlet problem in domains with {A}hlfors regular boundary}, 
      author={Ariel Barton and Svitlana Mayboroda and Alberto Pacati},
      year={2025},
      eprint={2506.14639},
      archivePrefix={arXiv},
      primaryClass={math.AP},
      url={https://arxiv.org/abs/2506.14639}, 
}

@Article{var77,
  author	= {Varopoulos, N. Th.},
  title		= {B{MO} functions and the {$\overline \partial $}-equation},
  journal	= {Pacific J. Math.},
  fjournal	= {Pacific Journal of Mathematics},
  volume	= {71},
  year		= {1977},
  number	= {1},
  pages		= {221--273},
  issn		= {0030-8730},
  mrclass	= {32A99 (46E99)},
  mrnumber	= {508035},
  mrreviewer	= {John E. Forn\ae ss},
  url		= {http://projecteuclid.org/euclid.pjm/1102811647}
}

@Article{var78,
  author	= {Varopoulos, N. Th.},
  title		= {A remark on functions of bounded mean oscillation and
		  bounded harmonic functions. {A}ddendum to: ``{BMO}
		  functions and the {$\overline \partial $}-equation''
		  ({P}acific {J}. {M}ath. {\bf 71} (1977), no. 1, 221--273)},
  journal	= {Pacific J. Math.},
  fjournal	= {Pacific Journal of Mathematics},
  volume	= {74},
  year		= {1978},
  number	= {1},
  pages		= {257--259},
  issn		= {0030-8730},
  mrclass	= {32A99 (46E99)},
  mrnumber	= {508036},
  mrreviewer	= {John E. Forn\ae ss},
  url		= {http://projecteuclid.org/euclid.pjm/1102810455}
}

@unpublished{mp26,
      title={The {D}irichlet problem as the boundary of the {P}oisson problem: {A} sharp approximation result}, 
      author={Mihalis Mourgoglou and Bruno Poggi}, 
      note={Preprint. February 2026. arXiv: 2602.07560},
      eprint={2602.07560},
      archivePrefix={arXiv},
      primaryClass={math.AP},
      url={https://arxiv.org/abs/2602.07560}, 
}

@article {hmt10,
    AUTHOR = {Hofmann, Steve and Mitrea, Marius and Taylor, Michael},
     TITLE = {Singular integrals and elliptic boundary problems on regular
              {S}emmes-{K}enig-{T}oro domains},
   JOURNAL = {Int. Math. Res. Not. IMRN},
  FJOURNAL = {International Mathematics Research Notices. IMRN},
      YEAR = {2010},
    NUMBER = {14},
     PAGES = {2567--2865},
      ISSN = {1073-7928,1687-0247},
   MRCLASS = {42B25 (35J25 42B20)},
  MRNUMBER = {2669659},
MRREVIEWER = {Brett\ D.\ Wick},
       DOI = {10.1093/imrn/rnp214},
       URL = {https://doi.org/10.1093/imrn/rnp214},
}

@misc{mt26,
      title={Solvability of the {N}eumann problem for elliptic equations in chord-arc domains with very big pieces of good superdomains}, 
      author={Mihalis Mourgoglou and Xavier Tolsa},
      year={2026},
      eprint={2407.20385},
      archivePrefix={arXiv},
      primaryClass={math.AP},
      url={https://arxiv.org/abs/2407.20385}, 
}

@article{jk81,
author = {David S. Jerison and Carlos E. Kenig},
title = {{The Neumann problem on Lipschitz domains}},
volume = {4},
journal = {Bulletin (New Series) of the American Mathematical Society},
number = {2},
publisher = {American Mathematical Society},
pages = {203 -- 207},
year = {1981},
}

@article {am14,
    AUTHOR = {Auscher, Pascal and Mourgoglou, Mihalis},
     TITLE = {Boundary layers, {R}ellich estimates and extrapolation of
              solvability for elliptic systems},
   JOURNAL = {Proc. Lond. Math. Soc. (3)},
  FJOURNAL = {Proceedings of the London Mathematical Society. Third Series},
    VOLUME = {109},
      YEAR = {2014},
    NUMBER = {2},
     PAGES = {446--482},
      ISSN = {0024-6115,1460-244X},
   MRCLASS = {35J57 (35A01 35B45 42B25 42B35)},
  MRNUMBER = {3254931},
MRREVIEWER = {Juha\ K.\ Kinnunen},
       DOI = {10.1112/plms/pdu008},
       URL = {https://doi.org/10.1112/plms/pdu008},
}

@article {mt24,
    AUTHOR = {Mourgoglou, Mihalis and Tolsa, Xavier},
     TITLE = {The regularity problem for the {L}aplace equation in rough
              domains},
   JOURNAL = {Duke Math. J.},
  FJOURNAL = {Duke Mathematical Journal},
    VOLUME = {173},
      YEAR = {2024},
    NUMBER = {9},
     PAGES = {1731--1837},
      ISSN = {0012-7094,1547-7398},
   MRCLASS = {42B37 (31B10 31B15 35J25 42B25 49Q15)},
  MRNUMBER = {4765702},
MRREVIEWER = {Juha\ K.\ Kinnunen},
       DOI = {10.1215/00127094-2023-0044},
       URL = {https://doi.org/10.1215/00127094-2023-0044},
}

@book {eg1,
    AUTHOR = {Evans, Lawrence C. and Gariepy, Ronald F.},
     TITLE = {Measure theory and fine properties of functions},
    SERIES = {Studies in Advanced Mathematics},
   EDITION = {1},
 PUBLISHER = {CRC Press, Boca Raton, FL},
      YEAR = {1992},
     PAGES = {xiv+268},
      ISBN = {0-8493-7157-0},
}

@article {kp,
    AUTHOR = {Kenig, Carlos E. and Pipher, Jill},
     TITLE = {The {N}eumann problem for elliptic equations with nonsmooth
              coefficients},
   JOURNAL = {Invent. Math.},
  FJOURNAL = {Inventiones Mathematicae},
    VOLUME = {113},
      YEAR = {1993},
    NUMBER = {3},
     PAGES = {447--509},
      ISSN = {0020-9910},
   MRCLASS = {35J25},
  MRNUMBER = {1231834},
MRREVIEWER = {Bruno Franchi},
       DOI = {10.1007/BF01244315},
       URL = {https://doi.org/10.1007/BF01244315},
}

@article {m,
    AUTHOR = {Mayboroda, Svitlana},
     TITLE = {The connections between {D}irichlet, regularity and {N}eumann
              problems for second order elliptic operators with complex
              bounded measurable coefficients},
   JOURNAL = {Adv. Math.},
  FJOURNAL = {Advances in Mathematics},
    VOLUME = {225},
      YEAR = {2010},
    NUMBER = {4},
     PAGES = {1786--1819},
      ISSN = {0001-8708},
   MRCLASS = {35B65 (35J25)},
  MRNUMBER = {2680190},
MRREVIEWER = {Tatiana Savina},
       DOI = {10.1016/j.aim.2010.04.019},
       URL = {https://doi.org/10.1016/j.aim.2010.04.019},
}

@article {mpt,
    AUTHOR = {Milakis, Emmanouil and Pipher, Jill and Toro, Tatiana},
     TITLE = {Harmonic analysis on chord arc domains},
   JOURNAL = {J. Geom. Anal.},
  FJOURNAL = {Journal of Geometric Analysis},
    VOLUME = {23},
      YEAR = {2013},
    NUMBER = {4},
     PAGES = {2091--2157},
      ISSN = {1050-6926},
   MRCLASS = {35B20 (31B35 46E30)},
  MRNUMBER = {3107693},
MRREVIEWER = {Srinivasan Natesan},
       DOI = {10.1007/s12220-012-9322-4},
       URL = {https://doi-org.ezp1.lib.umn.edu/10.1007/s12220-012-9322-4},
}

@article {a,
    AUTHOR = {Auscher, Pascal},
     TITLE = {Regularity theorems and heat kernel for elliptic operators},
   JOURNAL = {J. London Math. Soc. (2)},
  FJOURNAL = {Journal of the London Mathematical Society. Second Series},
    VOLUME = {54},
      YEAR = {1996},
    NUMBER = {2},
     PAGES = {284--296},
      ISSN = {0024-6107},
   MRCLASS = {35J15 (35B45 35K10 58G11)},
  MRNUMBER = {1405056},
MRREVIEWER = {Yehuda Pinchover},
       DOI = {10.1112/jlms/54.2.284},
       URL = {https://doi.org/10.1112/jlms/54.2.284},
}

@article {as,
    AUTHOR = {Aizenman, M. and Simon, B.},
     TITLE = {Brownian motion and {H}arnack inequality for {S}chr\"{o}dinger
              operators},
   JOURNAL = {Comm. Pure Appl. Math.},
  FJOURNAL = {Communications on Pure and Applied Mathematics},
    VOLUME = {35},
      YEAR = {1982},
    NUMBER = {2},
     PAGES = {209--273},
      ISSN = {0010-3640},
   MRCLASS = {35J10 (60J65 81C05)},
  MRNUMBER = {644024},
MRREVIEWER = {M. Schechter},
       DOI = {10.1002/cpa.3160350206},
       URL = {https://doi.org/10.1002/cpa.3160350206},
}

@unpublished {dm,
    AUTHOR = {David, Guy and Mayboroda, Svitlana},
     TITLE = {Harmonic measure is absolutely continuous with respect to the Hausdorff measure on all low-dimensional uniformly rectifiable sets},
      NOTE = {Preprint. arXiv:2006.14661},
      URL = {https://arxiv.org/abs/2006.14661},
}

@unpublished {dhp22,
    AUTHOR = {Dindo\v{s}, Martin and Hofmann, Steve and Pipher, Jill},
     TITLE = {Regularity and {N}eumann problems for operators with real coefficients satisfying {C}arleson condition},
      NOTE = {Preprint. July 2022 (v1), December 2022 (v2). arXiv: 2207.10366},
      URL = {https://arxiv.org/abs/2207.10366},
}

@article {mpt25,
    AUTHOR = {Mourgoglou, Mihalis and Poggi, Bruno and Tolsa, Xavier},
     TITLE = {Solvability of the {P}oisson-{D}irichlet problem with interior data in ${L}^{p'}$-{C}arleson spaces and its applications to the ${L}^p$-regularity problem},
   JOURNAL = {J. Eur. Math. Soc.},
  FJOURNAL = {Journal of the European Mathematical Society},
    VOLUME = {Published online first},
      YEAR = {2025},
       DOI = {10.4171/JEMS/1660},
       URL = {https://ems.press/journals/jems/articles/14298929},
}

@article{HS25,
  author  = {Hofmann, Steve and Sparrius, Derek},
  title   = {The Neumann function and the {$L^p$} Neumann problem in chord-arc domains},
  journal = {Advanced Nonlinear Studies},
  volume  = {25},
  number  = {2},
  pages   = {313--349},
  year    = {2025},
  doi     = {10.1515/ans-2023-0171}
}

@article{ABHM19,
  author  = {Akman, Murat and Bortz, Simon and Hofmann, Steve and Martell, Jos{\'e} Mar{\'i}a},
  title   = {Rectifiability, interior approximation and harmonic measure},
  journal = {Ark. Mat.},
  volume  = {57},
  year    = {2019},
  number  = {1},
  pages   = {1--22},
  doi     = {10.4310/ARKIV.2019.v57.n1.a1},
}

@article{M21,
  author  = {Mourgoglou, Mihalis},
  title   = {Approximate tangents, harmonic measure, and domains with rectifiable boundaries},
  journal = {Proc. Amer. Math. Soc.},
  volume  = {149},
  year    = {2021},
  number  = {9},
  pages   = {3739--3749},
  doi     = {10.1090/proc/15542},
}

\end{document}